\documentclass[12pt]{amsart}
\usepackage{placeins}
\usepackage{amssymb}
\newtheorem{theorem}{Theorem}
\newtheorem{example}[theorem]{Example}
\newtheorem{proposition}[theorem]{Proposition}

\newtheorem{corollary}[theorem]{Corollary}
\newtheorem{remark}[theorem]{Remark}
\newtheorem{definition}{Definition}
\newtheorem{lemma}[theorem]{Lemma}

\newcommand{\G}{{\mathfrak G}}

\newcommand{\HH}{{\mathfrak H}}
\newcommand{\Z}{{\mathbb Z}}

\newcommand{\R}{{\mathbb R}}
\newcommand{\C}{{\mathbb C}}

\begin{document}

\title[Caustics of $J_{10}$ singularities]{Complements of caustics of the real $J_{10}$ singularities}
\author{V.A.~Vassiliev}
\address{Weizmann Institute of Science, Rehovot, Israel}
\thanks{This work was supported by the Absorption Center in Science of the Ministry of Immigration and Absorption of the State of Israel}
 \email{vavassiliev@gmail.com}

\begin{abstract}
The complete list of connected components of the set of Morse functions in the deformations of function singularities of class $J_{10}$ is given. 
Thus, the isotopy classification of Morse perturbations of parabolic real function singularities is finished.
\end{abstract}

\maketitle

{

\section{Introduction}

\subsection{}
This work completes the isotopy classification of Morse perturbations of real parabolic function singularities. The geometry of sets of Morse perturbations of the
simplest singularity classes $A_2$, $A_3$, $A_4$ and $D_4^{\pm}$ was studied by R.~Thom and V.~Arnold in the context of catastrophe theoretical problems occurring in biology and optics, see \cite{thom}, \cite{thom2}, \cite{Acongr}, \cite{Kluwer}. 
The local components of the space of Morse functions near all
{\em simple} function singularities were enumerated in \cite{sed}, \cite{sede}, and \cite{Vsing}. For parabolic singularities of classes $X_9$ and $P_8$, the analogous problem was solved in \cite{Vx9} and \cite{Vp8}, respectively. Here, we solve the same problem for the remaining class of parabolic singularities, $J_{10}$. 
We prove that there are exactly 59 and 56 isotopy classes of Morse perturbations of $J_{10}^1$ and $J_{10}^3$ singularities, respectively.

In \S~\ref{sect2}, we recall a combinatorial invariant of isotopy classes of Morse functions (i.e., connected components of the space of these functions) that was used in previous works \cite{Vsing}--\cite{Vp8}. Then, we compute all possible values of this invariant for $J_{10} $ singularities and prove that each value can be realized by exactly one or two isotopy classes, depending on whether or not these classes are invariant or not under a certain symmetry. Finally, we examine this symmetry condition. Additionally, we realize many isotopy classes by concrete polynomials.

Our main invariant of isotopy classes of Morse functions is formulated in the terms of a graph, whose vertices correspond to the collections of certain topological characteristics of Morse functions, and whose edges correspond to their standard surgeries. The invariant's values are the subgraphs of this graph, into which it splits when the edges corresponding to those surgeries that can change the Morse isotopy class are removed. Also, the main
criterion for the self-symmetry condition of isotopy classes can be expressed in the terms of one-dimensional cocycles of this graph.

The most essential part of the obtained classification table is the list of the isotopy classes of polynomials with the maximal possible number of real critical points (ten), and only two distinct Morse indices. 
As with the $X_9$ and $P_8$ singularities, nearly all of these classes are related to the splittings of the original function singularity into pairs of real critical points, the sum of whose Milnor numbers is ten. The only additional quadruple of isotopy classes (all of which are mapped to each other by certain symmetries of the function space) is related in the same way to the canonical extended Coxeter-Dynkin diagram of class $\tilde E_8$ (another name for $J_{10}$ singularities). This situation repeats that of cases $X_9$ and $P_8$, where the unique exceptional isotopy classes are related to the diagrams of types $\tilde E_7$ and $\tilde E_6$.

\subsection{Main objects and definitions (see, e.g.,  \cite{AVGZ})}

A point $a \in \R^n$ is a {\em critical point} of a smooth function $\R^n \to \R$ if all first partial derivatives of this function vanish at $a$. 
A {\em function singularity} is a germ of a $C^\infty$-smooth function $(\R^n, a) \to (\R,0)$ at a point $a \in \R^n$ where its differential vanishes. 
 Two function singularities at points $a$ and $b$ are {\em equivalent} if they can be transformed to each other via a local diffeomorphism $(\R^n, a) \to (\R^n,b)$ (i.e., they have the same expression in appropriate local coordinate systems centered at $a$ and $b$). The {\em equivalence class} of a critical point $(f,a)$ is the equivalence class of the function singularity $(f-f(a),a).$ 
An $l$-parametric {\em deformation} of a function singularity $f$ is a family $\{f_\lambda \}$  of analytic functions depending on the parameter $\lambda$ running through a neighborhood of the origin point $0 \in \R^l$, such that $f_0 \equiv f$ and the function $F(x, \lambda) \equiv f_\lambda(x)$ in $n+l$ variables $x \in \R^n$ and $ \lambda \in \R^l$ is regular analytic. The {\em caustic} variety of the deformation $F \equiv \{f_\lambda\}$ is the set of parameter values $\lambda \in \R^l$ such that the corresponding function $f_\lambda$ has a non-Morse critical point near the origin in $\R^n$.

A smooth function $f:\R^2 \to \R$ has a $J_{10}$ singularity at a point $a \in \R^2$ if in some local coordinates $x$ and $ y$ centered at this point it is quasihomogeneous of degree six with weights $\deg x =2$ and $ \deg y=1$, and its Milnor number is finite (and then necessarily equal to ten, see \cite{AVGZ}, \cite{AGLV2}). The zero-level set $f^{-1}(0)$ of such a singularity in $\R^2$ can consist of one or three smooth local branches at \ $a$. \ The corresponding subclasses of the $J_{10}$ class are denoted by $J_{10}^1$ and $J_{10}^3$. 

\begin{table}
\caption{Real simple and $J_{10}$ singularities in two variables}
\label{t1}
\begin{center}
\begin{tabular}{|l|l|l|}
\hline
Notation & Normal form & Restriction \\ 
\hline 
$A_{2k-1}$ & $\pm x^{2k} \pm y^2 $ & $k \ge 1$ \\ 
$A_{2k}$ & $x^{2k+1} \pm y^2$ & $k \ge 1$ \\
[3pt]
$D_{k}^{\pm}$ & $x^2y \/\pm y^{k-1} $ & $k \ge 4$ \\ [4pt]
$E_6$ & $x^3 \pm y^4 $ & \cr
$E_7$ & $x^3 + x y^3 $ & \cr
$E_8$ & $x^3 + y^5$ & \cr
\hline
$J_{10}^1$ & $(x^2+y^4)(x-\gamma y^2)$ & $\gamma \in (-\infty, \infty)$ \\
$J_{10}^3$ & $(x^2-y^4)(x-\gamma y^2)$ & $ \gamma \in (-1,1)$ \\
\hline \end{tabular} \end{center}
\end{table}

The normal forms of $J_{10}$ singularities, to which they can be reduced by a choice of local coordinates, are shown in the last two rows of Table \ref{t1}. Other rows of this table show the normal forms of {\em simple singularities} into which the $J_{10}$ singularities can be moved by small perturbations.

The canonical {\em versal deformations} (i.e., ``sufficiently large'' deformations, to which all other deformations can be reduced, see \cite{AVGZ}) of function singularities of classes $J_{10}^1$ and $J_{10}^3$ consist respectively of polynomials
\begin{multline}
(x^2+y^4)(x-\gamma y^2) + \lambda_1 + \lambda_2 y + \lambda_3 y^2 + \lambda_4 y^3 + \lambda_5 y^4 + \\
\lambda_6 x +
\lambda_7 x y + \lambda_8 x y^2 + \lambda_9 (x^2+3y^4)y , \qquad \gamma \in (-\infty, +\infty) , \label{vers1} \end{multline}
and
\begin{multline}
(x^2-y^4)(x-\gamma y^2) + \lambda_1 + \lambda_2 y + \lambda_3 y^2 + \lambda_4 y^3 + \lambda_5 y^4 + \\
\lambda_6 x + \lambda_7 x y + \lambda_8 x y^2 + \lambda_9 (x^2-3y^4)y , \qquad \gamma \in (-1,1), \label{vers3}
\end{multline} 
 with ten parameters $\lambda_1, \dots, \lambda_9,$ and $\gamma$, see \cite{Jaw2}. 

The caustics in the parameter spaces $\R^9 \times \R$ and $\R^9 \times (-1,1)$ of these deformations divide them into connected components called {\em isotopy classes of Morse perturbations} of $J_{10}$ singularities. Our main goal is the enumeration of these components. 

As in \cite{Vx9}, \cite{Vp8}, we primarily consider slightly greater spaces than the canonical versal deformations (\ref{vers1}) and (\ref{vers3}), which can be reduced to these deformations by an appropriate group of diffeomorphisms $\R^2 \to \R^2$. 

Namely, we consider the 16-dimensional space of linear combinations of monomials $x^\alpha y^\beta$, such that $2 \alpha + \beta \leq 6$, with a positive coefficient on the monomial $x^3$ and a non-degenerate {\em principal quasihomogeneous part} (consisting of monomials $x^\alpha y^\beta$ with $2\alpha+\beta=6$). This principal part can be of class $J_{10}^1$ or $J_{10}^3$. The spaces of all polynomials with these principal parts are denoted, respectively, $\Phi_1$ and $\Phi_3$. Again, the sets of Morse functions of types $\Phi_1$ and $\Phi_3$ split into connected components (= isotopy classes). 

\begin{definition} \rm
The group $\G$ consists of all diffeomorphisms $\R^2 \to \R^2$ of the form 
\begin{equation}
\label{dig}
\tilde x = a x + b y+ c y^2 + \xi, \ \ \tilde y = d y + \eta, 
\end{equation}
 where $a>0, d>0. $ 
\end{definition}

This group is diffeomorphic to $\R^6$. It acts on the spaces $\Phi_1$ and $\Phi_3$ and preserves the set of Morse functions.  

\begin{proposition}
\label{proponne}
Each orbit of the action of the group $\G$ on the space $\Phi_1$ 
$($respectively, $\Phi_3)$ intersects the space of all polynomials of class
$($\ref{vers1}$)$ $($respectively, $($\ref{vers3}$))$
 transversally at a single point.
\end{proposition}

\noindent
{\it Proof.} 
It is easy to see that any polynomial of class $\Phi_1$ or $\Phi_3$ can be reduced to the form (\ref{vers1}) or (\ref{vers3}) by the action of the group $\G$. 
The Lie algebra of this group is generated by the vector fields $x \frac{\partial}{\partial x}$, $y \frac{\partial}{\partial x}$, $y^2 \frac{\partial}{\partial x}$, $\frac{\partial}{\partial x}$, $y \frac{\partial}{\partial y},$ and $\frac{\partial}{\partial y}$ corresponding to the infinitesimal changes of six parameters of this group. The Lie differentials of polynomials of class (\ref{vers1}) or  (\ref{vers3}) along these vector fields generate the normal bundles of these classes in the spaces $\Phi_1$ and $\Phi_3$. This implies the transversality of orbits to these classes.
Let $f$ be a polynomial of the form (\ref{vers1}), and let $G$ be an element of the group $\G$. Suppose the polynomial $\tilde f \equiv f \circ G$ also has the form (\ref{vers1}). Then, the  coefficients $a$ and $d$ of the diffeomorphism $G$ are equal to 1 and the coefficient $ c$ is equal to $0$, since otherwise the principal quasihomogeneous part of $\tilde f$ does not have the standard form of Table \ref{t1}. $\eta = 0$ since the coefficient of $\tilde f$ at the monomial $x y^2$ is trivial. $b = 0$ since otherwise the coefficient at $y^5$ is not equal to thrice the coefficient at $x^2 y$. $\xi = 0$ since the coefficient at $x^2$ is trivial. Thus, $G = \mbox{Id}$. The proof for deformation (\ref{vers3}) is analogous. \hfill $\Box$ \medskip

 So, each space $\Phi_1$ or $\Phi_3$ is canonically diffeomorphic to the direct product of the group $\G$ and the space of polynomials (\ref{vers1}) or (\ref{vers3}).
Additionally, there is a one-to-one correspondence between the isotopy classes of Morse functions in spaces $\Phi_1$ and $\Phi_3$ and the connected components of their intersections with the spaces of polynomials (\ref{vers1}) and (\ref{vers3}), respectively. 
\medskip

Alternatively, any function of type $\Phi_1$ or $\Phi_3$
can be reduced by an element of the group $\G$
 to the normal form with zero coefficients at the monomials $x^2, x^2y, $ and $x^2y^2$, coefficient 1 at $x^3$, and either coefficient $\pm 1$ at $y^6$ and zero coefficient at $y^5$, or zero coefficient at $y^6$ and $x y^3$ and coefficient $\pm 1 $ at $x y^4$. \medskip

The main result of the paper is a list of all isotopy classes of Morse polynomials of classes $\Phi_1$ and $\Phi_3$ or, equivalently, of the form (\ref{vers1}) and (\ref{vers3}). It is formulated in Theorems \ref{enu1}, \ref{enu3}, and \ref{cher}.
\medskip

As a byproduct, we list all possible splittings of real $J_{10}$ singularities into pairs of critical points, such that the sum of their Milnor numbers is ten.

\begin{definition}[cf. \cite{Jaw2}, \S 2] \rm
For any pair of simple singularity classes $\Xi$ and $\tilde \Xi$ with $\mu(\Xi) + \mu(\tilde \Xi) = 10$, the notation $\{\Xi + \tilde \Xi\}\rightsquigarrow J_{10}^1$ or 
$\{\Xi + \tilde \Xi\}\rightsquigarrow J_{10}^3$ 
means that there exists a smooth function $f:{\mathbb R}^2 \to {\mathbb R}$ with a singularity of class $J_{10}^1$ or $J_{10}^3$ at the origin and a one-parametric deformation \ $\Theta: {\mathbb R}^2 \times [0, \varepsilon) \to {\mathbb R}$, \ $\Theta (\cdot, 0) \equiv f,$ \ of this function, such that, for any $\tau \in (0, \varepsilon) $, the corresponding function $f_\tau \equiv \Theta (\cdot, \tau)$ has a critical point of class $\Xi$ and a critical point of class $\tilde \Xi$ in such a way that these two critical points depend continuously on $\tau$ and tend to the origin in $\R^2$ as $\tau$ tends to 0.
\end{definition}

\begin{table}
\caption{Approximation of $J_{10}$ singularities by bisingularities}
\label{tabadj}
\begin{tabular}{|c||l|l|}
\hline
Type & $J_{10}^1$ & $J_{10}^3$ \\
\hline
$E_8 + A_2$ & Yes & No \\
$D_8^+ + A_2$ & No & No \\
$D_8^- + A_2$ & No & Yes \\
$A_8 + A_2$ & Yes & No \\
$E_7 + A_3$ & No & Yes \\
$D^{\pm}_7 + A_3$ & Yes & No \\
$A_7 + A_3$ & Yes & Yes \\
\hline
\end{tabular} \quad
\begin{tabular}{|c||l|l|}
\hline
Type & $J_{10}^1$ & $J_{10}^3$ \\
\hline
$ E_6 + D_4^+$ & Yes & No \\
$E_6 + D_4^-$ & No & No \\
$E_6 + A_4$ & Yes & No \\
$D_6^\pm + D_4^\pm$ & No & No \\
$D_6^+ + A_4$ & Yes & No \\
$D_6^- + A_4$ & No & Yes \\
\phantom{$D_6$} & & \\
\hline
\end{tabular} \quad
\begin{tabular}{|c||l|l|}
\hline
Type & $J_{10}^1$ & $J_{10}^3$ \\
\hline
$A_6 + D_4^+$ & No & No \\
$A_6 + D_4^- $ & No & Yes \\
$A_6 + A_4$ & Yes & No \\
$A_5 + A_5$ & Yes & Yes \\
$A_5 + D_5^{\pm}$ & Yes & Yes \\
$D_5^{\pm} + D_5^{\pm}$ & No & Yes \\
\phantom{$D_6$} & & \\
\hline
\end{tabular}
\end{table}

\begin{theorem}[see \S~\ref{sectlast}]
\label{tabadjp}
For any pair of simple real singularity classes $\Xi$ and $\tilde \Xi$ with $\mu(\Xi) + \mu(\tilde \Xi) = 10$, we have $\{\Xi + \tilde \Xi\}\rightsquigarrow J_{10}^1$
$($respectively, $\{\Xi + \tilde \Xi\}\rightsquigarrow J_{10}^3)$
 if and only if ``Yes'' is written in Table~\ref{tabadj} at the intersection of the row $\Xi + \tilde \Xi$ and the column $J_{10}^1$ $($respectively, $J_{10}^3)$. 
\end{theorem}

This theorem is a ``real'' analog of the corresponding part of Table III of \cite{Jaw2}.

\subsection{Lyashko--Looijenga map}
\label{LLM}
 The {\em Lyashko--Looijenga map} is
one of main tools of the proofs and computations in this work.  It maps a complexification of the spaces of polynomials (\ref{vers1}) and (\ref{vers3}) to the space $\mbox{Sym}^{10}(\C^1),$ taking each polynomial to the unordered collection of its critical values. This complexified space consists of all polynomials of the form (\ref{vers3}) with arbitrary complex coefficients and unique condition $\gamma \neq \pm1$. In the restriction to the set of complex polynomials with ten different critical values, this map defines a covering over the configuration space $B(\C^1, 10)$ of all cardinality 10 subsets in $\C^1$, see \cite{Jaw}, \cite{Jaw2}. 

Both real spaces (\ref{vers1}) and (\ref{vers3}) can be embedded into this complex space. The embedding of the space (\ref{vers3}) is literal.
In the case of the space of real polynomials (\ref{vers1}), we first identify its parameter space with a real subspace of this complex deformation using the substitution $y = e^{\pi i/4} \tilde y$. Namely, this subspace consists of polynomials of the form  (\ref{vers3}) with coefficients $(\lambda, \gamma)$ such that the numbers $i \gamma,$ $ \lambda_1,$ $ e^{\pi i/4} \lambda_2,$ $ i \lambda_3,$ $ e^{3 \pi i/4} \lambda_4,$ $ \lambda_5, $ $\lambda_6, e^{\pi i/4}\lambda_7,$ $ i \lambda_8,$ and $ \lambda_9$ are real. These polynomials are considered as real functions of the real subspace in $\C^2$ consisting of points $(x,y)$ where $x \in \R$, $e^{-\pi i/4} y \in \R$.

The restrictions of the Lyashko--Looijenga map to these real spaces send the real polynomials (\ref{vers1}) and (\ref{vers3}) to the space of point collections in $\C^1$ that are invariant under the complex conjugation.  Over the collections, all of whose points are distinct, these maps are local diffeomorphisms. Also, according to \cite{Jaw2}, they behave in a proper way near the simplest degenerations of polynomials. This allows us to reduce the topology of  spaces of Morse polynomials to the combinatorics of point configurations and associated structures.

\section{Invariants of isotopy classes of Morse polynomials}
\label{sect2}

This section is a shortened version of Section 2 of \cite{Vx9}.

\subsection{Trivial invariant} 
\label{itriv}

The simplest invariant of isotopy classes of Morse polynomials in ${\mathbb R}^2$ is their {\it passport}, i.e., the triple $(m_-, m_{\times}, m_+)$ of the numbers of their minima, saddlepoints, and maxima. 

According to index considerations, for any Morse polynomial of type $\Phi_1$ or $\Phi_3$ the sum $M \equiv m_- + m_{\times} + m_+$ is an even number no greater than ten, and the {\em Euler number} $m_- - m_{\times} + m_+$ is equal to 0 for all Morse polynomials of type $\Phi_1$, and to $-2$ for polynomials of type $\Phi_3$. Therefore, when studying polynomials of a particular type $\Phi_1$ or $\Phi_3$, we will express the passports by only pairs of numbers, $m_+$ and $M$. These numbers determine the remaining passport numbers.

\subsection{Set-valued invariant and virtual Morse functions}

\label{svinv}

\begin{definition} \rm
A polynomial $f: ({\mathbb C}^2, {\mathbb R}^2) \to ({\mathbb C}, {\mathbb R})$ of type $\Phi_1$ or $\Phi_3$ is {\em generic} if it has only Morse critical points in ${\mathbb C}^2$, and all the corresponding critical values are different and not equal to 0.
\end{definition}

We associate a set of discrete topological characteristics, called a {\em virtual Morse function},
with any generic Morse polynomial $f$. 

\unitlength 1.00mm
\linethickness{0.4pt}
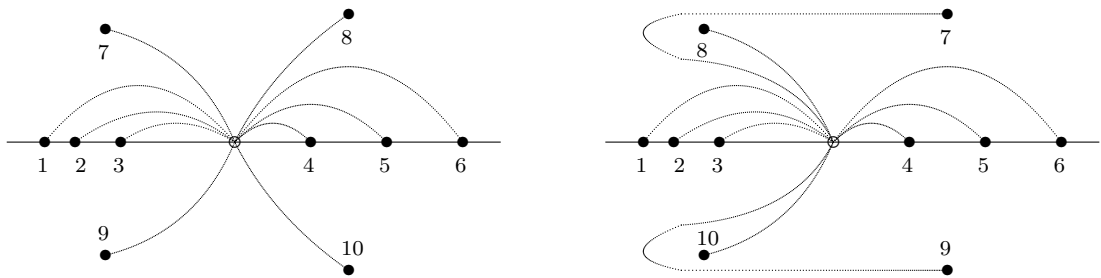
\begin{figure}
\begin{center}
\begin{picture}(65,40)
\put(0,20){\line(1,0){65}}
\put(5,20){\circle*{1.5}}
\put(9,20){\circle*{1.5}}
\put(15,20){\circle*{1.5}}
\put(40,20){\circle*{1.5}}
\put(50,20){\circle*{1.5}}
\put(60,20){\circle*{1.5}}
\bezier{100}(30,20)(17,35)(5,20)
\bezier{80}(30,20)(20,28)(9,20)
\bezier{50}(30,20)(22,25)(15,20)
\bezier{50}(30,20)(35,25)(40,20)
\bezier{90}(30,20)(40,30)(50,20)
\bezier{140}(30,20)(45,40)(60,20)
\put(13,35){\circle*{1.5}}
\put(13,5){\circle*{1.5}}
\put(45,37){\circle*{1.5}}
\put(45,3){\circle*{1.5}}
\bezier{100}(30,20)(25,32)(13,35)
\bezier{100}(30,20)(25,8)(13,5)
\bezier{100}(30,20)(35,30)(45,37)
\bezier{100}(30,20)(35,10)(45,3)
\put(30,20){\circle{1.5}}
\put(4,16){{\tiny 1}}
\put(9,16){{\tiny 2}}
\put(14,16){{\tiny 3}}
\put(39,16){{\tiny 4}}
\put(49,16){{\tiny 5}}
\put(59,16){{\tiny 6}}
\put(12,31){{\tiny 7}}
\put(44,33){{\tiny 8}}
\put(12,7){{\tiny 9}}
\put(44,5){{\tiny 10}}
\end{picture} \qquad \quad
\begin{picture}(65,40)
\put(0,20){\line(1,0){65}}
\put(5,20){\circle*{1.5}}
\put(9,20){\circle*{1.5}}
\put(15,20){\circle*{1.5}}
\put(40,20){\circle*{1.5}}
\put(50,20){\circle*{1.5}}
\put(60,20){\circle*{1.5}}
\bezier{100}(30,20)(17,35)(5,20)
\bezier{80}(30,20)(20,28)(9,20)
\bezier{50}(30,20)(22,25)(15,20)
\bezier{50}(30,20)(35,25)(40,20)
\bezier{90}(30,20)(40,30)(50,20)
\bezier{140}(30,20)(45,40)(60,20)
\put(13,35){\circle*{1.5}}
\put(13,5){\circle*{1.5}}
\put(45,37){\circle*{1.5}}
\put(45,3){\circle*{1.5}}
\bezier{100}(30,20)(25,32)(13,35)
\bezier{100}(30,20)(25,8)(13,5)
\bezier{100}(30,20)(25,30)(10,31)
\bezier{50}(10,31)(0,35)(10,37)
\bezier{100}(10,37)(35,37)(45,37)
\bezier{100}(30,20)(25,10)(10,9)
\bezier{50}(10,9)(0,5)(10,3)
\bezier{100}(10,3)(35,3)(45,3)
\put(30,20){\circle{1.5}}
\put(4,16){{\tiny 1}}
\put(9,16){{\tiny 2}}
\put(14,16){{\tiny 3}}
\put(39,16){{\tiny 4}}
\put(49,16){{\tiny 5}}
\put(59,16){{\tiny 6}}
\put(12,31.3){{\tiny 8}}
\put(44,33){{\tiny 7}}
\put(12,6.4){{\tiny 10}}
\put(44,5){{\tiny 9}}
\end{picture}
\end{center}
\caption{Standard systems of paths}
\label{standd}
\end{figure}

If $f$ is a generic polynomial of type $\Phi_1$ or $\Phi_3$, then the set $V_f \subset {\mathbb C}^3$ defined by the equation $f(x, y) + z^2=0$ is a smooth complex surface that is homotopy equivalent to the wedge of ten two-dimensional spheres (see, for example, \cite{M} and \cite{AVGZ}). The homology group $H_2(V_f)$ is generated by {\em vanishing cycles} (see \cite{AVGZ}, \cite{APLT}) that are defined by a system of non-intersecting paths in ${\mathbb C}^1$ connecting the non-critical value 0 with all critical values of $f$, see Fig.~\ref{standd}. We choose these paths so that those going to complex conjugate non-real critical values are symmetric about the real axis, and those going to real values lie in the upper half-plane where the imaginary parts are smaller than the absolute values of the imaginary parts of all non-real critical values. 

Let us fix an orientation of ${\mathbb R}^3$ somehow. Then there is a canonical choice of orientations of all vanishing cycles defined by such a system of paths, see \S V.1.6 of \cite{APLT}. In particular, complex conjugation in ${\mathbb C}^3$ must take the oriented vanishing cycles defined by complex conjugate paths into each other with the coefficient 1 and not $-1$.

 An order of these vanishing cycles can also be canonically defined. In particular, the cycles that vanish at real critical points are listed first in ascending order of the corresponding critical values.

\begin{definition}[see \cite{AGLV2}, \S V.3]\rm
\label{vm1}

A {\em virtual Morse function associated with} a generic Morse {\em polynomial $f:({\mathbb C}^2, {\mathbb R}^2) \to ({\mathbb C}, {\mathbb R})$} of type $\Phi_1$ or $\Phi_3$ is a collection of its topological data consisting of

a) the $10 \times 10$ matrix of intersection indices in $V_f$
of canonically ordered and oriented vanishing cycles $\Delta_i \in H_2(V_f)$ corresponding to all critical values of $f$ and defined by a system of paths as above,

b) the string of ten intersection indices in $V_f$ of these vanishing cycles with the naturally oriented set of real points, 
$V_f \cap {\mathbb R}^3$;

c) the string of {\em positive} Morse indices (i.e., positive inertia indices of the quadratic parts) of all real critical points of the function $f(x,y)+z^2$, and

d) the numbers of negative and positive real critical values of $f$.
\end{definition}

\noindent
{\bf Example.} 
Some four virtual Morse functions with eight real critical points are shown in Fig.~\ref{cocy}. Two vertical lines in each table indicate the last element of the corresponding virtual Morse function. Namely, in all four cases, the numbers of negative, positive and non-real critical values are 3, 5, and 2, respectively. The real critical points of these functions are only minima and saddlepoints, see the bottom lines of the tables.

{
\normalsize
\begin{figure}
\begin{equation*}
\begin{array}{|ccc|ccccc|cc|}
\hline
 $-2$ & $0$ & $0$ & $0$ & $0$ & $1$ & $0$ & $0$ & $0$ & $0$ \\
 $0$ & $-2$ & $0$ & $0$ & $1$ & $0$ & $1$ & $0$ & $0$ & $0$ \\
 $0$ & $0$ & $-2$ & $1$ & $0$ & $1$ & $0$ & $1$ & $1$ & $-1$ \\
 $0$ & $0$ & $1$ & $-2$ & $0$ & $0$ & $0$ & $0$ & $-1$ & $0$ \\
 $0$ & $1$ & $0$ & $0$ & $-2$ & $0$ & $0$ & $0$ & $-1$ & $-1$ \\
 $1$ & $0$ & $1$ & $0$ & $0$ & $-2$ & $0$ & $0$ & $-1$ & $0$ \\
 $0$ & $1$ & $0$ & $0$ & $0$ & $0$ & $-2$ & $0$ & $0$ & $0$ \\
 $0$ & $0$ & $1$ & $0$ & $0$ & $0$ & $0$ & $-2$ & $0$ & $1$ \\
 $0$ & $0$ & $1$ & $-1$ & $-1$ & $-1$ & $0$ & $0$ & $-2$ & $-1$ \\ 
 $0$ & $0$ & $-1$ & $0$ & $-1$ & $0$ & $0$ & $1$ & $-1$ & $-2$ \\
\hline 
 $0$ & $0$ & $0$ & $-1$ & $-1$ & $0$ & $-1$ & $-1$ & $-1$ & $-1$ \\
\hline
 $3$ & $3$ & $3$ & $2$ & $2$ & $2$ & $2$ & $2$ & & \\
\hline
\end{array} \ \ \ \
\begin{array}{|ccc|ccccc|cc|}
\hline
 $-2$ & $0$ & $0$ & $0$ & $0$ & $1$ & $0$ & $0$ & $0$ & $0$ \\
 $0$ & $-2$ & $0$ & $0$ & $1$ & $1$ & $1$ & $0$ & $-1$ & $1$ \\
 $0$ & $0$ & $-2$ & $1$ & $0$ & $0$ & $0$ & $1$ & $0$ & $0$ \\
 $0$ & $0$ & $1$ & $-2$ & $0$ & $0$ & $0$ & $0$ & $1$ & $1$ \\
 $0$ & $1$ & $0$ & $0$ & $-2$ & $0$ & $0$ & $0$ & $1$ & $0$ \\
 $1$ & $1$ & $0$ & $0$ & $0$ & $-2$ & $0$ & $0$ & $1$ & $0$ \\
 $0$ & $1$ & $0$ & $0$ & $0$ & $0$ & $-2$ & $0$ & $0$ & $-1$ \\
 $0$ & $0$ & $1$ & $0$ & $0$ & $0$ & $0$ & $-2$ & $0$ & $0$ \\
 $0$ & $-1$ & $0$ & $1$ & $1$ & $1$ & $0$ & $0$ & $-2$ & $-1$ \\
 $0$ & $1$ & $0$ & $1$ & $0$ & $0$ & $-1$ & $0$ & $-1$ & $-2$ \\
\hline
 $0$ & $0$ & $0$ & $-1$ & $-1$ & $0$ & $-1$ & $-1$ & $1$ & $1$ \\
\hline
 $3$ & $3$ & $3$ & $2$ & $2$ & $2$ & $2$ & $2$ & & \\
\hline
\end{array}
\end{equation*} 

\begin{equation*}
\begin{array}{|ccc|ccccc|cc|}
\hline 
 $-2$ & $0$ & $0$ & $0$ & $0$ & $1$ & $0$ & $0$ & $0$ & $0$ \\
 $0$ & $-2$ & $0$ & $1$ & $0$ & $1$ & $0$ & $1$ & $1$ & $-1$ \\
 $0$ & $0$ & $-2$ & $0$ & $1$ & $0$ & $1$ & $0$ & $0$ & $0$ \\
 $0$ & $1$ & $0$ & $-2$ & $0$ & $0$ & $0$ & $0$ & $-1$ & $0$ \\
 $0$ & $0$ & $1$ & $0$ & $-2$ & $0$ & $0$ & $0$ & $-1$ & $-1$ \\
 $1$ & $1$ & $0$ & $0$ & $0$ & $-2$ & $0$ & $0$ & $-1$ & $0$ \\
 $0$ & $0$ & $1$ & $0$ & $0$ & $0$ & $-2$ & $0$ & $0$ & $0$ \\
 $0$ & $1$ & $0$ & $0$ & $0$ & $0$ & $0$ & $-2$ & $0$ & $1$ \\
 $0$ & $1$ & $0$ & $-1$ & $-1$ & $-1$ & $0$ & $0$ & $-2$ & $-1$ \\
 $0$ & $-1$ & $0$ & $0$ & $-1$ & $0$ & $0$ & $1$ & $-1$ & $-2$ \\
\hline 
 $0$ & $0$ & $0$ & $-1$ & $-1$ & $0$ & $-1$ & $-1$ & $-1$ & $-1$ \\
\hline
 $3$ & $3$ & $3$ & $2$ & $2$ & $2$ & $2$ & $2$ & & \\
\hline 
\end{array} \ \ \ \
\begin{array}{|ccc|ccccc|cc|}
\hline
 $-2$ & $0$ & $0$ & $0$ & $0$ & $1$ & $0$ & $0$ & $0$ & $0$ \\
 $0$ & $-2$ & $0$ & $0$ & $1$ & $1$ & $0$ & $1$ & $1$ & $-1$ \\
 $0$ & $0$ & $-2$ & $1$ & $0$ & $0$ & $1$ & $0$ & $0$ & $0$ \\
 $0$ & $0$ & $1$ & $-2$ & $0$ & $0$ & $0$ & $0$ & $-1$ & $-1$ \\
 $0$ & $1$ & $0$ & $0$ & $-2$ & $0$ & $0$ & $0$ & $-1$ & $0$ \\
 $1$ & $1$ & $0$ & $0$ & $0$ & $-2$ & $0$ & $0$ & $-1$ & $0$ \\
 $0$ & $0$ & $1$ & $0$ & $0$ & $0$ & $-2$ & $0$ & $0$ & $0$ \\
 $0$ & $1$ & $0$ & $0$ & $0$ & $0$ & $0$ & $-2$ & $0$ & $1$ \\
 $0$ & $1$ & $0$ & $-1$ & $-1$ & $-1$ & $0$ & $0$ & $-2$ & $-1$ \\ 
 $0$ & $-1$ & $0$ & $-1$ & $0$ & $0$ & $0$ & $1$ & $-1$ & $-2$ \\
\hline 
 $0$ & $0$ & $0$ & $-1$ & $-1$ & $0$ & $-1$ & $-1$ & $-1$ & $-1$ \\
\hline
 $3$ & $3$ & $3$ & $2$ & $2$ & $2$ & $2$ & $2$ & & \\
\hline
\end{array}
\end{equation*}
\caption{Virtual Morse functions}
\label{cocy}
\end{figure}
}

\begin{definition} \rm 
A {\em critical point of a virtual Morse function} is any column of its data set as in Fig.~\ref{cocy}, i.e., a column of the intersection matrix, the intersection index with the set of real points, and a Morse index or the information that the critical point is non-real.
\end{definition}

\begin{remark} \rm
If a real polynomial $f$ has more than one pair of non-real critical values, then there can be more than one virtual Morse function associated with $f$, because the choice of a proper system of paths is not homotopically unique: see Fig.~\ref{standd}.
\end{remark}

\begin{definition} \rm
\label{elsur}
{\em Elementary virtual surgeries}\index{elementary surgery} of virtual Morse functions include six transformations of their data, modeling the standard local topological surgeries of the corresponding generic Morse polynomials, namely
\begin{itemize}
\item[$s1, s2$:] \ \ \ collision of two neighboring real critical values at a non-zero value, after which the corresponding two critical points either ($s1$) meet and leave the real domain, or ($s2$) change the order in ${\mathbb R}^1$ of their critical values; 

\item[$s3, s4$:] \ \ \ collision of two complex conjugate critical values at a point on the line ${\mathbb R}^1 \setminus \{0\}$, after which the corresponding critical points either ($s3$) meet at a real point and enter real space, or ($s4$) miss each other in the complex domain, while the imaginary parts of their critical values change their signs; 

\item[$s5, s6$:] \ \ \ jumps of real critical values up ($s5$) or down ($s6$) through 0; \\ \hspace*{4cm} and additionally

\item[$s7$:] \ \ \ specifically virtual transformations within the classes of virtual Morse functions associated with the same real Morse polynomials, that are caused by flips of standard systems of paths going from 0 to non-real critical values (see Fig.~\ref{standd} and also Figs. 19--21 of \cite{AVGZ}, volume 2).
\end{itemize}
\end{definition}

The results of all these virtual surgeries are determined by the data of the initial virtual Morse functions. For a detailed description of these standard flips of data, see \S V.8 of \cite{APLT}. The explicit formulas for them are given in the comments of the computer program that performs them, see the link on page \pageref{progg}.
In particular, attempting to perform the surgery $s1$ or $s2$ over real critical values $v_i$ and $v_{i+1}$ begins with examining the intersection index $\langle \Delta_i, \Delta_{i+1}\rangle$ of the corresponding vanishing cycles. If this 
index is $0$, then surgery $s2$ occurs; if the index is 1, then surgery $s1$ occurs; in all other cases the surgery fails. Similarly, a collision of two complex conjugate critical values at a real point not separated from 0 by other critical values follows scenario $s4$ if the intersection index is 0, scenario $s3$ if the index is $1$ or $-1$, and fails in all other cases; in the second case the sign of the intersection index allows us to predict the Morse indices of the newborn real critical points.

We will denote by $s1, \dots, s6$ both the real surgeries of real Morse functions and the corresponding elementary virtual surgeries.

\begin{remark} \rm
If our deformed singularities were neither parabolic nor simple,  we could not be sure that the real surgery corresponding to a virtual one could be realized at any time when these intersection index conditions are satisfied. For simple singularities this realization is guaranteed by the properness of the Lyashko--Looijenga map (see \cite{Lo0}). For parabolic singularities, this follows from the P.~Jaworski's work \cite{Jaw}, \cite{Jaw2}, see Proposition \ref{propmain} below.
\end{remark}

Let $f: {\mathbb R}^2 \to {\mathbb R}$ be a generic Morse polynomial of type $\Phi_1$ or $\Phi_3$.

\begin{definition} \rm
An (abstract) {\em virtual Morse function} of type $f$ is any collection of data as in Definition \ref{vm1} (i.e., a matrix, two strings, and two numbers) obtained from an arbitrary virtual Morse function associated with $f$ via an arbitrary finite chain of elementary virtual surgeries. 

The {\em formal graph} of type $f$ is the graph, whose vertices correspond to all virtual Morse functions of type $f$, and two such vertices are connected by an edge if and only if the corresponding virtual Morse functions can be obtained from each other by an elementary virtual surgery.

The {\em virtual component} $S(f)$ of the formal graph of type $f$ is its subgraph, whose vertices are only the virtual Morse functions of type $f$ that can be obtained from virtual Morse functions associated with $f$ via arbitrary finite chains of virtual surgeries $s2, s4, s5, s6,$ and $s7$ from Definition \ref{elsur} (i.e., all surgeries that do not model the collision of critical points).
\end{definition}

\begin{example} \rm
\label{cocycleh}
Denote the four virtual Morse functions in Fig.~\ref{cocy} by the letters A, B, C, and D in accordance with the diagram 
${A \ B \atop C \ D}$. The transition $A \leftrightarrow B$ describes the passage $s4$ of two imaginary critical values through the real axis between the fifth and sixth real critical values. 
The transitions $A \leftrightarrow C$, $B \leftrightarrow D$ and $C\leftrightarrow D$ are  the surgeries $s2$ at which respectively
the second and third, the seventh and eighth, and the fourth and fifth real critical values meet and overtake each other. These four elementary surgeries form a commutative diagram, and the chain of passages $A \to B \to D \to C \to A$ forms a cycle in a virtual component of the formal graph. 
\end{example}

\begin{proposition}
\label{protrivi}
For all generic polynomials $f$ of the same type $\Phi_1$ or $\Phi_3$,
 the formal graph of type $f$ is the same.

If two generic Morse polynomials $f$ and $\tilde f$ belong to the same connected component of the set of Morse polynomials of type $\Phi_1$ or $\Phi_3$, then their virtual components $S(f)$ and $S(\tilde f)$ are the same.
\end{proposition}

\noindent
{\it Proof.} The first statement follows from the connectedness of the spaces $\Phi_1$ and $\Phi_3$. Any two real Morse polynomials from such a space can be connected by a path that intersects the variety of non-Morse polynomials at finitely many points, each of which corresponds to a standard surgery. The second statement follows immediately from the definitions. \hfill $\Box$
\medskip

Proposition \ref{protrivi} justifies the following definitions.

\begin{definition} \rm
\label{def8}
The {\em formal graph} of a (not necessarily Morse) polynomial ${\mathbb R}^2 \to {\mathbb R}$ of type $\Phi_1$ or $\Phi_3$ is the formal graph of an arbitrary generic Morse polynomial of the same type. 

The {\it virtual component} of a (not necessarily strictly) Morse polynomial of type $\Phi_1$ or $\Phi_3$ is the virtual component of an arbitrary generic polynomial from the same connected component of the space of all Morse polynomials of this type.
\end{definition}

\begin{proposition}
\label{propmain}
For any generic Morse polynomial $f$ of type $\Phi_1$ or $\Phi_3$, a virtual Morse function $\varphi$ associated with it, and a virtual Morse function $\tilde \varphi \neq \varphi$ connected with $\varphi$ by an edge of the formal graph not of type $s7$, there exists a generic Morse polynomial $\tilde f$ associated with $\tilde \varphi$ and a path in the space $\Phi_1$ or $\Phi_3$ connecting $f$ and $\tilde f$ and containing only one non-generic point at which it experiences a
standard surgery of the same type as the edge $[\varphi, \tilde \varphi]$.
\end{proposition}

\noindent
{\it Proof}. For type $\Phi_3$ polynomials, the proof repeats the proof of Proposition 2 in \cite{Vx9} for $X_9^+$ singularities. The only difference is that it considers the canonical  versal deformation of the complex $J_{10}$ singularities, which has exactly the form (\ref{vers3}) with complex parameters and the condition $\gamma \neq \pm 1$ instead of the canonical versal deformation of $X_9$. For type $\Phi_1$, we first embed 
the deformation (\ref{vers1}) into the complex $J_{10}$ deformation as described in \S \ref{LLM}. \hfill $\Box$

\begin{corollary}
\label{cormain}
Let $f$ be a generic polynomial of type $\Phi_1$ or $\Phi_3$, then

a$)$ every virtual Morse function of type $f$ is associated with a generic real polynomial $\tilde f$ of the same type $\Phi_1$ or $\Phi_3$ as $f$;

b$)$ every virtual Morse function from the virtual component $S(f)$ is associated with some generic real polynomial from the same connected component as $f$ of the space of Morse polynomials of type $\Phi_1$ or $\Phi_3$. \hfill $\Box$
\end{corollary}

\begin{definition} \rm
\label{svi}
The {\em set-valued invariant} of a real Morse polynomial $f$ of type $\Phi_1$ or $\Phi_3$ is the set of virtual Morse functions corresponding to all vertices of the virtual component $S(f)$. The {\em invariant} \ Card \ of such a Morse polynomial $f$ is the cardinality of this set of vertices.
\end{definition}

Clearly, the set-valued invariant determines the ``passport'' invariant (which can be read from the bottom line of any virtual Morse function).

\begin{proposition}
\label{uniq}
If two generic Morse polynomials of type $\Phi_1$ or $\Phi_3$ are associated with the same virtual Morse function, then they either belong to the same connected component of the space of generic polynomials of this type, or they belong to the connected components that are mapped to each other by the involution 
\begin{equation}
f(x,y) \leftrightarrow f(x,-y) . \label{invol0}
\end{equation}
\end{proposition}

The proof of this proposition will be given in \S~\ref{secproof}.

\begin{corollary}
Each virtual component of the formal graph of type $\Phi_1$ or $\Phi_3$ is associated with one or two isotopy classes of Morse polynomials of the corresponding type $\Phi_1$ or $\Phi_3$. \hfill $\Box$
\end{corollary}

\begin{definition} \rm
\label{her}
An isotopy class of Morse polynomials and the virtual component associated with this isotopy class are called {\em achiral} (respectively, {\em chiral}) if the involution (\ref{invol0}) takes this isotopy class to itself (respectively, to a different component). 
\end{definition}

Finally, our system of isotopy invariants of Morse polynomials of types $\Phi_1$ and $\Phi_3$ consists of the set-valued invariant of Definition \ref{svi} and the reflection class in the case of chiral classes. 
According to Proposition \ref{uniq}, this system of invariants separates all isotopy classes.

\begin{proposition}
\label{proinvar}
If a polynomial $f$ of type $\Phi_1$ or $\Phi_3$ is invariant under the involution $($\ref{invol0}$)$, then the intersection of its $\G$-orbit with the space of polynomials $($\ref{vers1}$)$ or $($\ref{vers3}$)$ is also invariant under this involution.
\end{proposition}

\noindent
{\it Proof.} This intersection point and its image under this involution belong to the same $\G$-orbit. According to Proposition \ref{proponne}, this orbit has only one intersection point with the space (\ref{vers1}) or (\ref{vers3}). \hfill $\Box$

\subsection{D-graph invariant}
\label{dinv}

In this subsection we only consider the real polynomials $({\mathbb C}^2,{\mathbb R}^2) \to ({\mathbb C}, {\mathbb R}),$ whose critical points are all real. In this case the set-valued invariant of \S \ref{svinv} has the following transparent interpretation.

Let $f$ be a generic polynomial ${\mathbb R}^2 \to {\mathbb R}$ of type $\Phi_1$ or $\Phi_3$ with only real critical points, in particular all ten of its critical values are real and distinct, and $0$ is a non-critical value of $f$. The matrix of intersection indices of vanishing cycles $\Delta_i \in H_2(V_f)$ (numbered in the ascending order of the corresponding critical values)
 can be depicted by its {\em Coxeter-Dynkin graph} (see e.g. \cite{AVGZ}) with ten ordered vertices corresponding to all critical values of $f$. Namely, if the intersection index $\langle \Delta_i, \Delta_j \rangle$ is positive, then the corresponding vertices $v_i$ and $v_j$ are connected by $ \langle \Delta_i, \Delta_j \rangle$ solid segments; if $\langle \Delta_i, \Delta_j \rangle$ is negative then they are connected by $-\langle \Delta_i, \Delta_j \rangle$ dashed segments. 

\begin{definition}[see \cite{Vx9}] \rm
\label{dfdinv}The {\it D-graph} of a generic real Morse polynomial $f$ with only real critical points is (the isomorphism class of) the 
oriented graph with vertices labeled with indices 0, 1 and 2, which is obtained from the Coxeter-Dynkin graph of $f$ via

1) orienting each edge of the graph from the vertex corresponding to the critical point with the lower critical value to the vertex corresponding to the critical point with the higher critical value;

2) labeling each vertex of the graph with the Morse index of the corresponding critical point of $f$, and

3) forgetting the ordering of the vertices.
\end{definition} 

\noindent{\bf Notation.}
In the pictures of the D-graphs (see Figs.~\ref{122298}--\ref{29370}), instead of numerical Morse indices, 
we will label the vertices corresponding to minima, saddlepoints and maxima by white circles, black circles, and white squares, respectively. 

\begin{remark} \rm
The D-graph of a real Morse polynomial is determined by an arbitrary virtual Morse function associated with this polynomial. Indeed, the intersection matrix and the Morse indices are the elements of the virtual Morse function, and the orientation of the edges follows from the order of the rows and columns of the intersection matrix, which is determined by the order of the corresponding critical values. In this way, D-graphs of arbitrary {\em virtual} Morse functions with only real critical points are also well-defined.
\end{remark}

\begin{theorem}
In the restriction to the space of generic Morse polynomials of class $\Phi_1$ or $\Phi_3$
with only real critical points, the D-graphs form an invariant of isotopy classes of Morse functions. This
 invariant is equivalent to the set-valued invariant from \S \ref{svinv}.
\end{theorem}

\noindent 
{\it Proof} of this theorem repeats the proof of Theorem 2 of \cite{Vx9}. \hfill $\Box$

\subsection{The up-down involution}

The involution \begin{equation} 
f(x,y) \leftrightarrow -f(-x,y) \label{invol} \ \end{equation}
acts on the spaces $\Phi_1$ and $\Phi_3$. 

This action can also be extended to the corresponding virtual Morse functions and 
D-graphs. Namely, for any system of paths for the function $f(x,y)$ (see Fig.~\ref{standd}) we take the system of paths for $-f(-x,y)$ obtained from it by the composition of the multiplication by $-1$ and the complex conjugation in $\C^1$. The virtual Morse function associated with the function $-f(-x,y)$ and defined by the resulting system of paths is determined by that of the original virtual Morse function associated with $f(x,y)$, see \S~2.2 in \cite{Vp8}.

In particular, if two functions with ten real critical points are related via this involution, then their $D$-graphs are obtained from each other by replacing all minima by maxima and vice versa, and reversing the directions of all edges, see Proposition 16 in \cite{Vp8}.

\section{Enumeration of isotopy classes and virtual components}

\begin{table}
\caption{Numbers of virtual and real components of class $\Phi_1$ (left) and $\Phi_3$ (right)}
\label{J101}
\begin{tabular}{|c| c | c | c | c | c |c |}
\hline
$m_+ \backslash M$ & 0 & 2 & 4 & 6 & 8 & 10 \\
\hline
0 & 1 & 1 & 1 & 1 & 1(2) & 9(18) \\
1 & 0 & 1 & 1 & 1 & 1 & 1(2) \\
2 & 0 & 0 & 1 & 1 & 1 & 1 \\
3 & 0 & 0 & 0 & 1 & 1 & 1 \\
4 & 0 & 0 & 0 & 0 & 1(2) & 1(2) \\
5 & 0 & 0 & 0 & 0 & 0 & 9(18) \\
\hline
$\Sigma$ & 1 & 2 & 3 & 4 & 5(7) & 22(42) \\
\hline
\end{tabular} \qquad 
\begin{tabular}{|c| c | c | c | c | c | c |}
\hline
$m_+ \backslash M$ & 0 & 2 & 4 & 6 & 8 & 10 \\
\hline
0 & 0 & 1 & 1 & 1 & 2 & 9(18) \\
1 & 0 & 0 & 1 & 1 & 1 & 2(3) \\
2 & 0 & 0 & 0 & 1 & 1 & 1(2) \\
3 & 0 & 0 & 0 & 0 & 2 & 2(3) \\
4 & 0 & 0 & 0 & 0 & 0 & 9(18) \\
\hline 
$\Sigma$ & 0 & 1 & 2 & 3 & 6 & 23 (44) \\
\hline
\end{tabular}
\end{table}

\subsection{}
Table \ref{J101} shows the number of virtual components and the number of isotopy classes of Morse polynomials of types $\Phi_1$ and $\Phi_3$ with any value of the passport $(M, m_+)$. The second number is given in parentheses when it differs from the first. This table is detailed in Theorems \ref{enu1}, \ref{enu3}, \ref{cher}.

\begin{figure}
\unitlength=0.45mm
\begin{picture}(120,70)
\bezier{600}(60,30)(30,50)(15,40)
\bezier{600}(60,30)(30,10)(15,20)
\bezier{400}(15,20)(0,30)(15,40)
\bezier{400}(105,20)(120,30)(105,40)
\bezier{600}(60,30)(90,50)(105,40)
\bezier{600}(60,30)(90,10)(105,20)
\bezier{600}(15,0)(20,60)(32.5,60)
\bezier{600}(50,0)(45,60)(32.5,60)
\put(10,28){\footnotesize $-$}
\put(29,28){\footnotesize $+$}
\put(29,49){\footnotesize $-$}
\put(49,28){\footnotesize $-$}
\put(85,28){\footnotesize $-$}
\end{picture} \qquad \quad
\unitlength 1.2mm
\begin{picture}(53,27)
\put(0,13){\circle{1}}
\put(16,13){{\makebox(0,0)[cc]{$\diamond$}}}
\put(32,13){\circle{1}}
\put(43,13){\circle*{1}}
\put(54,13){\circle{1}}
\put(8,5){\circle*{1}}
\put(16,29){\circle{1}}
\put(24,5){\circle*{1}}
\put(8,21){\circle*{1}}
\put(24,21){\circle*{1}}
\put(0.4,12.6){\vector(1,-1){7.1}}
\put(0.4,13.4){\vector(1,1){7.1}}
\put(3,13){\line(1,0){4}}
\put(9,13){\vector(1,0){4.5}}
\put(15.6,28.6){\vector(-1,-1){7.1}}
\put(16.4,28.6){\vector(1,-1){7.1}}
\put(16,26){\line(0,-1){4}}
\put(16,19){\vector(0,-1){4.5}}
\put(8.4,20.6){\vector(1,-1){7.1}}
\put(23.6,20.6){\vector(-1,-1){7.1}}
\put(31.6,12.6){\vector(-1,-1){7.1}}
\put(31.6,13.4){\vector(-1,1){7.1}}
\put(8.4,5.4){\vector(1,1){7.1}}
\put(23.6,5.4){\vector(-1,1){7.1}}
\put(29,13){\line(-1,0){4}}
\put(22,13){\vector(-1,0){4.5}}
\put(32.6,13){\vector(1,0){9.7}}
\put(53.4,13){\vector(-1,0){9.7}}
\end{picture}
\caption{$\Phi_1$, one maximum (122298)}
\label{122298}
\end{figure}
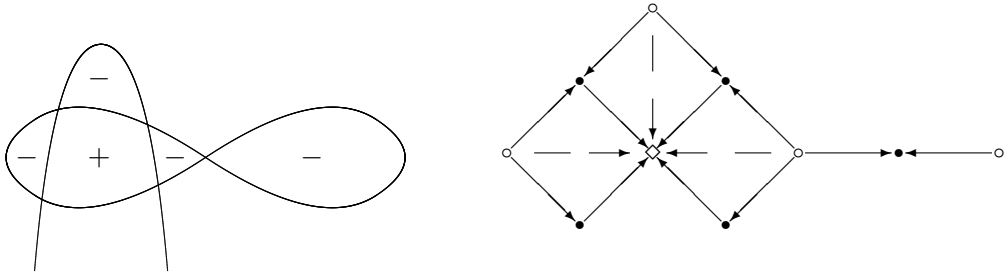

\begin{figure}
\unitlength 0.45mm
\begin{picture}(120,70)
\bezier{600}(60,30)(30,50)(15,40)
\bezier{600}(60,30)(30,10)(15,20)
\bezier{400}(15,20)(0,30)(15,40)
\bezier{400}(105,20)(120,30)(105,40)
\bezier{600}(60,30)(90,50)(105,40)
\bezier{600}(60,30)(90,10)(105,20)
\bezier{600}(30,0)(40,60)(60,60)
\bezier{600}(90,0)(80,60)(60,60)
\put(20,28){\footnotesize $-$}
\put(43,28){\footnotesize $+$}
\put(57,41){\footnotesize $-$}
\put(72,28){\footnotesize $+$}
\put(95,28){\footnotesize $-$}
\end{picture} \qquad \quad
\unitlength 0.45 mm
\begin{picture}(140,70)
\put(0,35){\circle{2}}
\put(40,35){{\makebox(0,0)[cc]{$\diamond$}}}
\put(70,35){\circle*{2}}
\put(100,35){{\makebox(0,0)[cc]{$\diamond$}}}
\put(140,35){\circle{2}}
\put(20,5){\circle*{2}}
\put(20,65){\circle*{2}}
\put(120,5){\circle*{2}}
\put(120,65){\circle*{2}}
\put(70,65){\circle{2}}
\put(0.5,35.7){\vector(2,3){18.7}}
\put(0.5,34.3){\vector(2,-3){18.7}}
\put(5,35){\line(1,0){8}}
\put(16,35){\line(1,0){8}}
\put(27,35){\vector(1,0){8}}
\put(20,5){\vector(2,3){19}}
\put(20,65){\vector(2,-3){19}}
\put(70,35){\vector(-1,0){28.5}}
\put(70,35){\vector(1,0){28.5}}
\put(139.5,35.7){\vector(-2,3){18.65}}
\put(139.5,34.3){\vector(-2,-3){18.65}}
\put(135,35){\line(-1,0){8}}
\put(124,35){\line(-1,0){8}}
\put(113,35){\vector(-1,0){8}}
\put(120,5){\vector(-2,3){19}}
\put(120,65){\vector(-2,-3){19}}
\put(70.7,65){\vector(1,0){48}}
\put(69.3,65){\vector(-1,0){48}}
\put(70,64.5){\vector(0,-1){29}}
\put(72,63){\line(1,-1){8}}
\put(82,53){\line(1,-1){8}}
\put(91,44){\vector(1,-1){8}}
\put(68,63){\line(-1,-1){8}}
\put(58,53){\line(-1,-1){8}}
\put(49,44){\vector(-1,-1){8}}
\end{picture}
\caption{$\Phi_1$, two maxima (26378)}
\label{26378}
\end{figure}
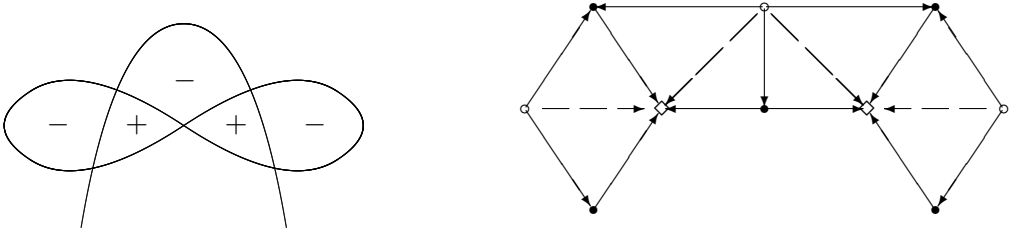

\begin{theorem}[$\Phi_1$]
\label{enu1} There are exactly 2005366 different virtual Morse functions of $J_{10}^1$ class distributed as follows into 37 virtual components.

\noindent 
1. There are exactly 22 virtual components of type $\Phi_1$ consisting of virtual Morse functions with ten real critical points:
\begin{itemize}
\item one component consisting of virtual Morse functions with exactly one local maximum, D-graph shown in Fig.~\ref{122298} $($right$)$ and Card invariant equal to 122298,
\item one component of virtual Morse functions with exactly two local maxima, D-graph shown in Fig.~\ref{26378} $($right$)$ and Card=26378,
\item nine components of virtual Morse functions without local maxima: their D-graphs are shown in Figs.~\ref{8499}--\ref{33528}, and their Card invariants are indicated in the captions of these figures, and
\item eleven virtual components obtained from the ones listed above in this theorem via 
the involution $($\ref{invol}$)$ of the space $\Phi_1$.
\end{itemize}
2. There are exactly five virtual components of type $\Phi_1$ with eight real critical points: 
\begin{itemize} 
\item one component with two local maxima, two minima, and Card=10890;
\item one component with exactly one local maximum and Card=27378, 
\item one component without local maxima, having Card=166554,
\item two virtual components obtained from the last two by the involution $($\ref{invol}$)$.
\end{itemize}
3. There are exactly four virtual components of type $\Phi_1$ with six real critical points: 
\begin{itemize} 
\item one component with exactly one local maximum and Card=10122, 
\item one component without local maxima, having Card=44128,
\item two virtual components obtained from the previous two by the involution $($\ref{invol}$)$.
\end{itemize}
4. There are exactly three virtual components of type $\Phi_1$ with four real critical points: 
\begin{itemize} 
\item one component with exactly one local maximum and Card=7850,
\item one component with no local maxima and Card=15850, 
\item a virtual component with no local minima obtained from the previous component via involution $($\ref{invol}$)$.
\end{itemize}
5. There are exactly two virtual components of type $\Phi_1$ with two real critical points: 
 one component without local maxima and one component without local minima obtained from it by involution $($\ref{invol}$)$, both having Card=10608.

\noindent
6. There is only one virtual component of class $\Phi_1$ without real critical points, its Card invariant is equal to 17642.
\end{theorem}

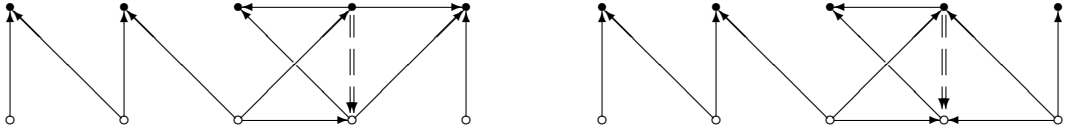
\begin{figure}
\unitlength 0.5 mm
\begin{picture}(120,50)
\put(0,5){\circle{2}}
\put(30,5){\circle{2}}
\put(60,5){\circle{2}}
\put(90,5){\circle{2}}
\put(120,5){\circle{2}}
\put(0,35){\circle*{2}}
\put(30,35){\circle*{2}}
\put(60,35){\circle*{2}}
\put(90,35){\circle*{2}}
\put(120,35){\circle*{2}}
\put(0,6){\vector(0,1){28}}
\put(29.3,5.7){\vector(-1,1){28.6}}
\put(30,6){\vector(0,1){28}}
\put(59.3,5.7){\vector(-1,1){28.6}}
\put(61,5){\vector(1,0){28}}
\put(60.7,5.7){\vector(1,1){28.6}}
\put(89,35){\vector(-1,0){28}}
\put(74.5,20.5){\vector(-1,1){13.8}}
\put(75.5, 19.5){\line(1,-1){13.8}}
\put(89.5,14){\vector(0,-1){7}}
\put(90.5,14){\vector(0,-1){7}}
\put(89.5,17){\line(0,1){7}}
\put(90.5,17){\line(0,1){7}}
\put(89.5,27){\line(0,1){6}}
\put(90.5,27){\line(0,1){6}}
\put(91,35){\vector(1,0){28}}
\put(90.7,5.7){\vector(1,1){28.6}}
\put(120,6){\vector(0,1){28}}
\end{picture} \qquad \qquad
\begin{picture}(120,50)
\put(0,5){\circle{2}}
\put(30,5){\circle{2}}
\put(60,5){\circle{2}}
\put(90,5){\circle{2}}
\put(120,5){\circle{2}}
\put(0,35){\circle*{2}}
\put(30,35){\circle*{2}}
\put(60,35){\circle*{2}}
\put(90,35){\circle*{2}}
\put(120,35){\circle*{2}}
\put(0,6){\vector(0,1){28}}
\put(29.3,5.7){\vector(-1,1){28.6}}
\put(30,6){\vector(0,1){28}}
\put(59.3,5.7){\vector(-1,1){28.6}}
\put(61,5){\vector(1,0){28}}
\put(60.7,5.7){\vector(1,1){28.6}}
\put(89,35){\vector(-1,0){28}}
\put(89.3,5.7){\line(-1,1){13.8}}
\put(74.5, 20.5){\vector(-1,1){13.7}}
\put(89.5,14){\vector(0,-1){6}}
\put(90.5,14){\vector(0,-1){6}}
\put(89.5,17){\line(0,1){7}}
\put(90.5,17){\line(0,1){7}}
\put(89.5,27){\line(0,1){6}}
\put(90.5,27){\line(0,1){6}}
\put(119,5){\vector(-1,0){28}}
\put(119.3,5.7){\vector(-1,1){28.6}}
\put(120,6){\vector(0,1){28}}
\end{picture}
\caption{$\Phi_1$, no maxima; 97702 (left) and 93489 (right)}
\label{8499}
\end{figure}

\begin{figure}
\unitlength 0.5 mm
\begin{picture}(120,50)
\put(0,5){\circle{2}}
\put(30,5){\circle{2}}
\put(60,5){\circle{2}}
\put(90,5){\circle{2}}
\put(120,5){\circle{2}}
\put(0,35){\circle*{2}}
\put(30,35){\circle*{2}}
\put(60,35){\circle*{2}}
\put(90,35){\circle*{2}}
\put(120,35){\circle*{2}}
\put(0,6){\vector(0,1){28}}
\put(29.3,5.7){\vector(-1,1){28.6}}
\put(30,6){\vector(0,1){28}}
\put(89.2,5.5){\vector(-2,1){58.4}}
\put(43,22){\vector(-1,1){9}}
\put(57,8){\line(-1,1){10}}
\put(120,6){\vector(0,1){28}}
\put(90.7,5.7){\vector(1,1){28.6}}
\put(91,35){\vector(1,0){28}}
\put(89.5,14){\vector(0,-1){7}}
\put(90.5,14){\vector(0,-1){7}}
\put(89.5,17){\line(0,1){7}}
\put(90.5,17){\line(0,1){7}}
\put(89.5,27){\line(0,1){6}}
\put(90.5,27){\line(0,1){6}}
\put(89,35){\vector(-1,0){28}}
\put(89.7,5.7){\vector(-1,1){28.6}}
\put(61,5){\vector(1,0){28}}
\bezier{500}(90,35)(50,45)(32.3,36.3)
\put(34,37){\vector(-2,-1){3}}
\put(75.6,20.6){\vector(1,1){13.8}}
\bezier{100}(74.6,19.6)(72.5,17.5)(70.5,15.5)
\put(69.5,14.5){\line(-1,-1){8.7}}
\end{picture} \qquad \qquad
\begin{picture}(120,50)
\put(0,5){\circle{2}}
\put(30,5){\circle{2}}
\put(60,5){\circle{2}}
\put(90,5){\circle{2}}
\put(120,5){\circle{2}}
\put(0,35){\circle*{2}}
\put(30,35){\circle*{2}}
\put(60,35){\circle*{2}}
\put(90,35){\circle*{2}}
\put(120,35){\circle*{2}}
\put(0,6){\vector(0,1){28}}
\put(29.3,5.7){\vector(-1,1){28.6}}
\put(30,6){\vector(0,1){28}}
\put(59.7,5.7){\vector(-1,1){28.6}}
\put(61,5){\vector(1,0){28}}
\put(60.7,5.7){\vector(1,1){28.6}}
\put(89,35){\vector(-1,0){28}}
\put(74.5,20.5){\vector(-1,1){13.8}}
\put(75.5, 19.5){\line(1,-1){13.8}}
\put(89.5,14){\vector(0,-1){7}}
\put(90.5,14){\vector(0,-1){7}}
\put(89.5,17){\line(0,1){7}}
\put(90.5,17){\line(0,1){7}}
\put(89.5,27){\line(0,1){6}}
\put(90.5,27){\line(0,1){6}}
\put(119,5){\vector(-1,0){28}}
\put(119.7,5.7){\vector(-1,1){28.6}}
\put(105.5,20.5){\vector(1,1){13.8}}
\put(104.5,19.5){\line(-1,-1){13.8}}
\put(91,35){\vector(1,0){28}}
\put(120,26){\vector(0,1){7}}
\put(120,23){\line(0,-1){6}}
\put(120,14){\line(0,-1){6.5}} 
\end{picture}
\caption{$\Phi_1$, no maxima; 68145 (left) and 42372 (right)}
\label{3852}
\end{figure}
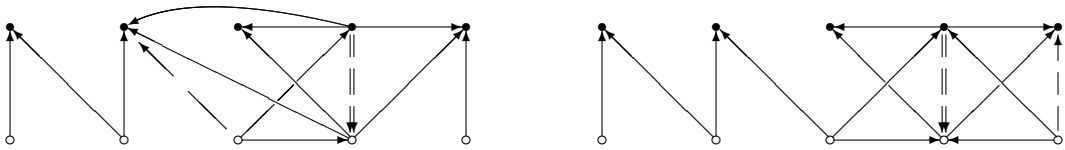

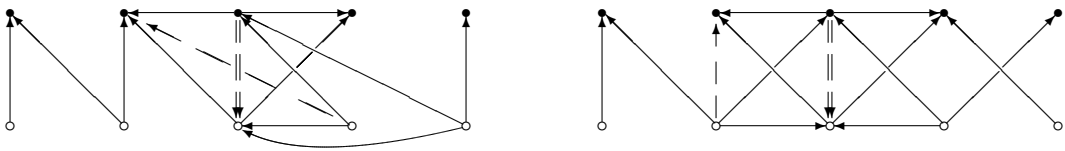
\begin{figure}
\unitlength 0.5 mm
\begin{picture}(120,50)
\put(0,5){\circle{2}}
\put(30,5){\circle{2}}
\put(60,5){\circle{2}}
\put(90,5){\circle{2}}
\put(120,5){\circle{2}}
\put(0,35){\circle*{2}}
\put(30,35){\circle*{2}}
\put(60,35){\circle*{2}}
\put(90,35){\circle*{2}}
\put(120,35){\circle*{2}}
\put(0,6){\vector(0,1){28}}
\put(29.3,5.7){\vector(-1,1){28.6}}
\put(30,6){\vector(0,1){28}}
\put(59.7,5.7){\vector(-1,1){28.6}}
\put(89,5){\vector(-1,0){28}}
\put(59.5,14){\vector(0,-1){7}}
\put(60.5,14){\vector(0,-1){7}}
\put(59.5,17){\line(0,1){7}}
\put(60.5,17){\line(0,1){7}}
\put(59.5,27){\line(0,1){6}}
\put(60.5,27){\line(0,1){6}}
\put(58,21){\line(-2,1){9}}
\put(45,27.5){\vector(-2,1){9}}
\put(62,19){\line(2,-1){8}}
\put(76,12){\line(2,-1){9}}


\put(59,35){\vector(-1,0){28}}
\put(61,35){\vector(1,0){28}}
\put(119.2,5.5){\vector(-2,1){58.4}}
\put(120,6){\vector(0,1){28}}
\put(89.3,5.3){\vector(-1,1){28.6}}
\put(80.5,25.5){\vector(1,1){8.8}}
\bezier{100}(79.5,24.5)(77.5,22.5)(75.5,20.5)
\put(74.3,19.3){\line(-1,-1){13.5}}
\bezier{300}(119,4.5)(80,-5)(63,3)
\put(63,3){\vector(-2,1){2}}
\end{picture} \qquad \qquad
\begin{picture}(120,50)
\put(0,5){\circle{2}}
\put(30,5){\circle{2}}
\put(60,5){\circle{2}}
\put(90,5){\circle{2}}
\put(120,5){\circle{2}}
\put(0,35){\circle*{2}}
\put(30,35){\circle*{2}}
\put(60,35){\circle*{2}}
\put(90,35){\circle*{2}}
\put(120,35){\circle*{2}}
\put(59.5,14){\vector(0,-1){7}}
\put(60.5,14){\vector(0,-1){7}}
\put(59.5,17){\line(0,1){7}}
\put(60.5,17){\line(0,1){7}}
\put(59.5,27){\line(0,1){6}}
\put(60.5,27){\line(0,1){6}}
\put(30,26){\vector(0,1){6}}
\put(30,22){\line(0,-1){6}}
\put(30,13){\line(0,-1){6}}
\put(0,6){\vector(0,1){28}}
\put(29.7,5.7){\vector(-1,1){28.7}}
\put(59,35){\vector(-1,0){28}}
\put(59.3,5.7){\vector(-1,1){28.6}}
\put(31,5){\vector(1,0){28}}
\put(89,5){\vector(-1,0){28}}
\put(89.7,5.7){\vector(-1,1){28.6}}
\put(61,35){\vector(1,0){28}}
\put(45.5,20.5){\vector(1,1){13.8}}
\put(44.5,19.5){\line(-1,-1){13.8}}
\put(75.5,20.5){\vector(1,1){13.8}}
\put(74.5,19.5){\line(-1,-1){13.8}}
\put(105.5,20.5){\vector(1,1){13.8}}
\put(104.5,19.5){\line(-1,-1){13.8}}
\put(119.3,5.7){\vector(-1,1){28.6}}
\end{picture}
\caption{$\Phi_1$, no maxima; 52415 (left) and 63085 (right)}
\label{5735}
\end{figure}

\begin{figure}
\unitlength 0.5 mm
\begin{picture}(120,50)
\put(0,5){\circle{2}}
\put(30,5){\circle{2}}
\put(60,5){\circle{2}}
\put(90,5){\circle{2}}
\put(120,5){\circle{2}}
\put(0,35){\circle*{2}}
\put(30,35){\circle*{2}}
\put(60,35){\circle*{2}}
\put(90,35){\circle*{2}}
\put(120,35){\circle*{2}}
\put(89.5,14){\vector(0,-1){7}}
\put(90.5,14){\vector(0,-1){7}}
\put(89.5,17){\line(0,1){7}}
\put(90.5,17){\line(0,1){7}}
\put(89.5,27){\line(0,1){6}}
\put(90.5,27){\line(0,1){6}}
\put(89,35){\vector(-1,0){28}}
\put(91,35){\vector(1,0){28}}
\put(120,6){\vector(0,1){28}}
\put(90.7,5.7){\vector(1,1){28.6}}
\put(43,22){\vector(-1,1){9}}
\put(57,8){\line(-1,1){10}}
\put(29.5,20.25){\vector(-2,1){28.4}}
\put(30.5,19.5){\line(2,-1){28.5}}
\put(30,6){\vector(0,1){28}}
\put(61,5){\vector(1,0){28}}
\put(89.3,5.7){\vector(-1,1){28.6}}
\put(89.2,5.5){\vector(-2,1){58.4}}
\bezier{400}(90,35)(50,45)(32,36)
\put(32,36){\vector(-2,-1){1}}
\put(0,6){\vector(0,1){28}}
\put(75.5,20.5){\vector(1,1){14}}
\bezier{100}(74.5,19.5)(72.5,17.5)(70.5,15.5)
\put(69.5,14.5){\line(-1,-1){8.8}}
\end{picture} \qquad \qquad
\begin{picture}(120,50)
\put(0,5){\circle{2}}
\put(30,5){\circle{2}}
\put(60,5){\circle{2}}
\put(90,5){\circle{2}}
\put(120,5){\circle{2}}
\put(0,35){\circle*{2}}
\put(30,35){\circle*{2}}
\put(60,35){\circle*{2}}
\put(90,35){\circle*{2}}
\put(120,35){\circle*{2}}
\put(89.5,14){\vector(0,-1){7}}
\put(90.5,14){\vector(0,-1){7}}
\put(89.5,17){\line(0,1){7}}
\put(90.5,17){\line(0,1){7}}
\put(89.5,27){\line(0,1){6}}
\put(90.5,27){\line(0,1){6}}
\put(43,22){\vector(-1,1){9}}
\put(57,8){\line(-1,1){10}}
\put(0,6){\vector(0,1){28}}
\put(29.3,5.7){\vector(-1,1){28.6}}
\put(30,6){\vector(0,1){28}}
\put(120,26){\vector(0,1){7}}
\put(120,22.5){\line(0,-1){6}}
\put(120,13){\line(0,-1){6}}
\put(119.3,5.7){\vector(-1,1){28.6}}
\put(61,5){\vector(1,0){28}}
\put(119,5){\vector(-1,0){28}}
\put(89,35){\vector(-1,0){28}}
\put(91,35){\vector(1,0){28}}
\put(89.3,5.7){\vector(-1,1){28.6}}
\put(89.2,5.5){\vector(-2,1){58.4}}
\put(105.5,20.5){\vector(1,1){13.8}}
\put(104.5,19.5){\line(-1,-1){14}}
\put(75.5,20.5){\vector(1,1){13.8}}
\bezier{100}(74.5,19.5)(72.5,17.5)(70.5,15.5)
\put(69.5, 14.5){\line(-1,-1){8.5}}
\bezier{400}(90,35)(50,45)(32,36)
\put(32,36){\vector(-2,-1){1}}
\end{picture}
\caption{$\Phi_1$, no maxima; 82500 (left) and 27940 (right)}
\label{2540}
\end{figure}
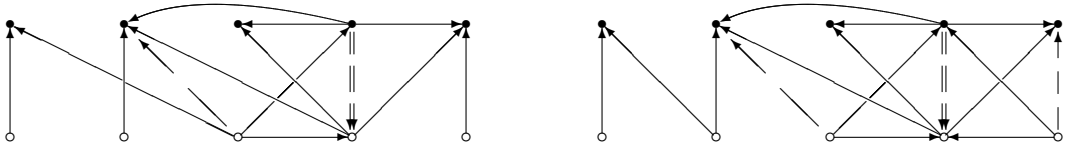

\begin{figure}
\unitlength 0.5 mm
\begin{picture}(120,50)
\put(0,5){\circle{2}}
\put(30,5){\circle{2}}
\put(60,5){\circle{2}}
\put(90,5){\circle{2}}
\put(120,5){\circle{2}}
\put(0,35){\circle*{2}}
\put(30,35){\circle*{2}}
\put(60,35){\circle*{2}}
\put(90,35){\circle*{2}}
\put(120,35){\circle*{2}}
\put(89.5,14){\vector(0,-1){7}}
\put(90.5,14){\vector(0,-1){7}}
\put(89.5,17){\line(0,1){7}}
\put(90.5,17){\line(0,1){7}}
\put(89.5,27){\line(0,1){6}}
\put(90.5,27){\line(0,1){6}}
\put(43,22){\vector(-1,1){9}}
\put(57,8){\line(-1,1){10}}
\put(0,6){\vector(0,1){28}}
\put(29.3,20.35){\vector(-2,1){28.2}}
\put(30.7,19.65){\line(2,-1){28.3}}
\put(30,6){\vector(0,1){28}}
\put(120,26){\vector(0,1){7}}
\put(120,22.5){\line(0,-1){6}}
\put(120,13){\line(0,-1){6}}
\put(119.3,5.7){\vector(-1,1){28.6}}
\put(61,5){\vector(1,0){28}}
\put(119,5){\vector(-1,0){28}}
\put(89,35){\vector(-1,0){28}}
\put(91,35){\vector(1,0){28}}
\put(89.3,5.7){\vector(-1,1){28.6}}
\put(89.2,5.4){\vector(-2,1){58.4}}
\put(105.5,20.5){\vector(1,1){13.8}}
\put(104.5,19.5){\line(-1,-1){14}}
\put(75.5,20.5){\vector(1,1){14.5}}
\bezier{100}(74.5,19.5)(72.5,17.5)(70.5,15.5)
\put(69.5, 14.5){\line(-1,-1){8.8}}
\bezier{400}(90,35)(50,45)(32,36)
\put(32,36){\vector(-2,-1){1}}
\end{picture}
\caption{$\Phi_1$, no maxima; 33528}
\label{33528}
\end{figure}

\FloatBarrier

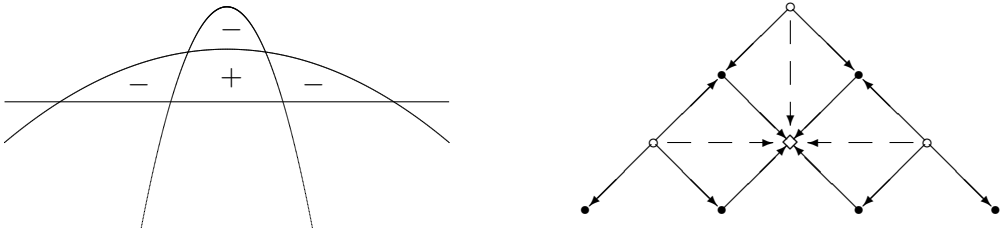
\begin{figure}
\unitlength 0.45mm
\begin{picture}(130,65)
\put(0,37){\line(1,0){130}}
\bezier{700}(0,25)(65,80)(130,25)
\bezier{600}(40,0)(65,130)(90,0)
\put(63,56){\footnotesize $-$}
\put(63,42){\footnotesize $+$}
\put(36,40){\footnotesize $-$}
\put(87,40){\footnotesize $-$}
\end{picture} \qquad \qquad \qquad
\begin{picture}(120,65)
\put(0,5){\circle*{2}}
\put(40,5){\circle*{2}}
\put(80,5){\circle*{2}}
\put(120,5){\circle*{2}}
\put(20,25){\circle{2}}
\put(60,25){{\makebox(0,0)[cc]{$\diamond$}}}
\put(100,25){\circle{2}}
\put(40,45){\circle*{2}}
\put(80,45){\circle*{2}}
\put(60,65){\circle{2}}
\put(19.3,24.3){\vector(-1,-1){18}}
\put(20.7,24.3){\vector(1,-1){18}}
\put(20.7,25.7){\vector(1,1){18}}
\put(40,5){\vector(1,1){18.5}}
\put(40,45){\vector(1,-1){18.5}}
\put(80,5){\vector(-1,1){18.5}}
\put(80,45){\vector(-1,-1){18.5}}
\put(99.3,24.3){\vector(-1,-1){18}}
\put(99.3,25.7){\vector(-1,1){18}}
\put(100.7,24.3){\vector(1,-1){18}}
\put(59.3,64.3){\vector(-1,-1){18}}
\put(60.7,64.3){\vector(1,-1){18}}
\put(24,25){\line(1,0){7}}
\put(36,25){\line(1,0){7}}
\put(48,25){\vector(1,0){7}}
\put(96,25){\line(-1,0){7}}
\put(84,25){\line(-1,0){7}}
\put(72,25){\vector(-1,0){7}}
\put(60,61){\line(0,-1){7}}
\put(60,49){\line(0,-1){7}}
\put(60,37){\vector(0,-1){7}}
\end{picture}
\caption{$\Phi_3$, one maximum (77374)}
\label{77374}
\end{figure}

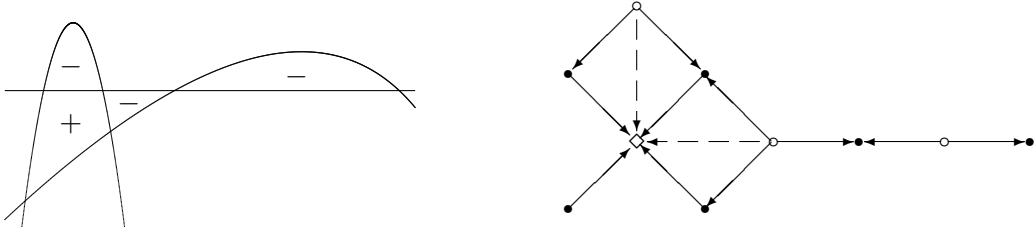
\begin{figure}
\unitlength 0.45mm
\begin{picture}(125,60)
\put(0,40){\line(1,0){120}}
\bezier{600}(5,0)(20,120)(35,0)
\bezier{800}(0,2)(83,80)(120,35)
\put(16,28){\footnotesize $+$}
\put(16,45){\footnotesize $-$}
\put(33,34){\footnotesize $-$}
\put(82,42){\footnotesize $-$}
\end{picture} \qquad \qquad \qquad
\begin{picture}(140,68)
\put(0,45){\circle*{2}}
\put(0,5){\circle*{2}}
\put(20,25){{\makebox(0,0)[cc]{$\diamond$}}}
\put(20,65){\circle{2}}
\put(40,45){\circle*{2}}
\put(40,5){\circle*{2}}
\put(60,25){\circle{2}}
\put(85,25){\circle*{2}}
\put(110,25){\circle{2}}
\put(135,25){\circle*{2}}
\put(0,5){\vector(1,1){18}}
\put(0,45){\vector(1,-1){18}}
\put(19.3,64.3){\vector(-1,-1){18}}
\put(20.7,64.3){\vector(1,-1){18}}
\put(20,61){\line(0,-1){6}}
\put(20,52){\line(0,-1){6}}
\put(20,43){\line(0,-1){6}}
\put(20,34){\vector(0,-1){6}}
\put(40,45){\vector(-1,-1){19}}
\put(40,5){\vector(-1,1){19}}
\put(59.3,25.3){\vector(-1,1){19}}
\put(59.3,24.7){\vector(-1,-1){19}}
\put(56,25){\line(-1,0){6}}
\put(47,25){\line(-1,0){6}}
\put(38,25){\line(-1,0){6}}
\put(29, 25){\vector(-1,0){6}}
\put(61,25){\vector(1,0){22.5}}
\put(109,25){\vector(-1,0){22.5}}
\put(111,25){\vector(1,0){22.5}}
\end{picture}
\caption{$\Phi_3$, one maximum (225148)}
\label{225148}
\end{figure}

\begin{figure}
\unitlength 0.45mm
\begin{picture}(125,65)
\put(0,45){\line(1,0){120}}
\bezier{600}(5,5)(40,135)(75,5)
\bezier{800}(0,12)(75,95)(120,30)
\put(19,38){\footnotesize $+$}
\put(33,55){\footnotesize $-$}
\put(50,47){\footnotesize $+$}
\put(73,49){\footnotesize $-$}
\end{picture} 
 \qquad \qquad \qquad
\begin{picture}(120,55)
\put(0,30){\circle*{2}}
\put(25,30){{\makebox(0,0)[cc]{$\diamond$}}}
\put(25,55){\circle*{2}}
\put(50,30){\circle*{2}}
\put(50,55){\circle{2}}
\put(75,30){{\makebox(0,0)[cc]{$\diamond$}}}
\put(100,30){\circle*{2}}
\put(75,55){\circle*{2}}
\put(100,55){\circle{2}}
\put(125,55){\circle*{2}}
\put(0,30){\vector(1,0){24}}
\put(25,55){\vector(0,-1){24}}
\put(47.3,52.3){\line(-1,-1){8}}
\put(36.3,41.3){\vector(-1,-1){8}}
\put(52.7,52.3){\line(1,-1){8}}
\put(63.7,41.3){\vector(1,-1){8}}
\put(50,30){\vector(-1,0){23}}
\put(50,54){\vector(0,-1){23}}
\put(49,55){\vector(-1,0){23}}
\put(51,55){\vector(1,0){23}}
\put(51,30){\vector(1,0){23}}
\put(75,55){\vector(0,-1){23}}
\put(100,30){\vector(-1,0){23}}
\put(101,55){\vector(1,0){23}}
\put(99,55){\vector(-1,0){23}}
\put(100,54){\vector(0,-1){23}}
\put(97.3,52.3){\line(-1,-1){8}}
\put(86.3,41.3){\vector(-1,-1){8}}
\end{picture}
\caption{$\Phi_3$, two maxima (128634)}
\label{128634}
\end{figure}
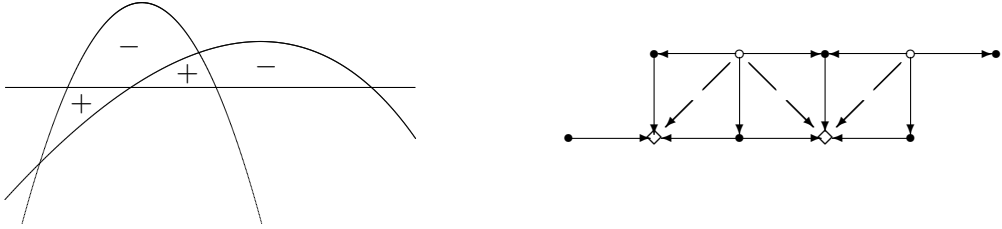

\begin{theorem}[$\Phi_3$]
\label{enu3}
There are exactly 2970134 different virtual Morse functions of $J_{10}^3$ type distributed as follows into 35 virtual components.

\noindent
1. There are exactly 23 virtual components of type $\Phi_3$ with ten real critical points:
\begin{itemize}
\item two components with exactly one local maximum: one with Card=77374 and D-graph shown in Fig.~\ref{77374} $($right$)$, and the other with Card=225148 and D-graph shown in Fig.~\ref{225148} $($right$)$;
\item one component with two local maxima, two minima, Card=128634 and D-graph shown in Fig.~\ref{128634} $($right$)$;
\item nine virtual components without local maxima, whose D-graphs are shown in Figs.~\ref{102234}--\ref{29370} and Card invariants are indicated in the captions of these pictures; 
\item eleven virtual components obtained by the involution $($\ref{invol}$)$ from all the components listed above in this theorem except for the component with two minima and two maxima.
\end{itemize}
2. There are exactly six virtual components of type $\Phi_3$ with eight real critical points: 
\begin{itemize} 
\item one component with one local maximum, two local minima, and Card=66906,
\item two components with no local maxima and values of Card invariant equal to 131148 and 82350,
\item three virtual components obtained from these mentioned in the previous two items by involution $($\ref{invol}$)$.
\end{itemize}
3. There are exactly three virtual components of type $\Phi_3$ with only six real critical points: 
\begin{itemize} 
\item one component with one local maximum, one local minimum, and Card=26922,
\item one component without local maxima having Card=57442, 
\item a virtual component obtained from the previous one by the involution $($\ref{invol}$)$.
\end{itemize}
4. There are exactly two virtual components of type $\Phi_3$ with only four real cri\-tical points: 
one component without local maxima, having Card=21410,
and one component obtained from it by the involution $($\ref{invol}$)$;

\noindent
5. There is only one virtual component of type $\Phi_3$ with exactly two real critical points: it has Card=14778 and no local extrema;

\noindent
6. There are no virtual components of type $\Phi_3$ without real critical points.
\end{theorem}

\begin{figure}
\unitlength 0.45mm
\begin{picture}(150,50)
\put(30,5){\circle{2}}
\put(60,5){\circle{2}}
\put(90,5){\circle{2}}
\put(120,5){\circle{2}}
\put(0,35){\circle*{2}}
\put(30,35){\circle*{2}}
\put(60,35){\circle*{2}}
\put(90,35){\circle*{2}}
\put(120,35){\circle*{2}}
\put(150,35){\circle*{2}}
\put(29.3,5.7){\vector(-1,1){28.6}}
\put(30,6){\vector(0,1){28}}
\put(59.3,5.7){\vector(-1,1){28.6}}
\put(60,6){\vector(0,1){28}}
\put(89.3,5.7){\vector(-1,1){28.6}}
\put(89,35){\vector(-1,0){28}}
\put(89.5,14){\vector(0,-1){7}}
\put(90.5,14){\vector(0,-1){7}}
\put(89.5,17){\line(0,1){7}}
\put(90.5,17){\line(0,1){7}}
\put(89.5,27){\line(0,1){6}}
\put(90.5,27){\line(0,1){6}}
\put(91,35){\vector(1,0){28}}
\put(90.7,5.7){\vector(1,1){28.6}}
\put(90.8,5.4){\vector(2,1){58.5}}
\put(120.7,5.7){\vector(1,1){28.6}}
\bezier{400}(150,35.5)(110,45)(91.5,36.5)
\put(91.5,36.5){\vector(-3,-2){1}}
\end{picture} \qquad 
\begin{picture}(150,50)
\put(30,5){\circle{2}}
\put(60,5){\circle{2}}
\put(90,5){\circle{2}}
\put(120,5){\circle{2}}
\put(0,35){\circle*{2}}
\put(30,35){\circle*{2}}
\put(60,35){\circle*{2}}
\put(90,35){\circle*{2}}
\put(120,35){\circle*{2}}\put(150,35){\circle*{2}}
\put(29.3,5.7){\vector(-1,1){28.6}}
\put(30,6){\vector(0,1){28}}
\put(59.3,5.7){\vector(-1,1){28.6}}
\put(60,6){\vector(0,1){28}}
\put(89.3,5.7){\vector(-1,1){28.6}}
\put(89,35){\vector(-1,0){28}}
\put(89.5,14){\vector(0,-1){7}}
\put(90.5,14){\vector(0,-1){7}}
\put(89.5,17){\line(0,1){7}}
\put(90.5,17){\line(0,1){7}}
\put(89.5,27){\line(0,1){6}}
\put(90.5,27){\line(0,1){6}}
\put(91,35){\vector(1,0){28}}
\put(90.7,5.7){\vector(1,1){28.6}}
\put(120.7,5.7){\vector(1,1){28.6}}
\put(119,5){\vector(-1,0){28}}
\put(104.5,20.5){\vector(-1,1){13.8}}
\put(105.5,19.5){\line(1,-1){13.9}}
\end{picture}
\caption{$\Phi_3$, no maxima; 230472 (left) and 102234 (right)}
\label{102234}
\end{figure}

\begin{figure}
\unitlength 0.45mm
\begin{picture}(150,50)
\put(30,5){\circle{2}}
\put(60,5){\circle{2}}
\put(90,5){\circle{2}}
\put(120,5){\circle{2}}
\put(0,35){\circle*{2}}
\put(30,35){\circle*{2}}
\put(60,35){\circle*{2}}
\put(90,35){\circle*{2}}
\put(120,35){\circle*{2}}
\put(150,35){\circle*{2}}
\put(59,35){\vector(-1,0){28}}
\put(61,35){\vector(1,0){28}}
\put(29.3,5.7){\vector(-1,1){28.6}}
\put(59.2,5.4){\vector(-2,1){58.4}}
\put(59.3,5.7){\vector(-1,1){28.6}}
\put(89.7,5.7){\vector(-1,1){28.6}}
\put(120.7,5.7){\vector(1,1){28.6}}
\put(119.3,5.7){\vector(-1,1){28.6}}
\bezier{400}(2,36)(20,45)(60,35)
\put(2,36){\vector(-2,-1){1}}
\put(75.5,20.5){\vector(1,1){13.8}}
\put(74.5,19.5){\line(-1,-1){14}}
\put(105.5,20.5){\vector(1,1){13.8}}
\put(104.5,19.5){\line(-1,-1){14}}
\put(90,26){\vector(0,1){6}}
\put(90,22){\line(0,-1){6}}
\put(90,12){\line(0,-1){6}}
\put(89,5){\vector(-1,0){28}}
\put(59.5,14){\vector(0,-1){7}}
\put(60.5,14){\vector(0,-1){7}}
\put(59.5,17){\line(0,1){7}}
\put(60.5,17){\line(0,1){7}}
\put(59.5,27){\line(0,1){6}}
\put(60.5,27){\line(0,1){6}}
\end{picture}
\qquad 
\begin{picture}(150,50)
\put(30,5){\circle{2}}
\put(60,5){\circle{2}}
\put(90,5){\circle{2}}
\put(120,5){\circle{2}}
\put(0,35){\circle*{2}}
\put(30,35){\circle*{2}}
\put(60,35){\circle*{2}}
\put(90,35){\circle*{2}}
\put(120,35){\circle*{2}}
\put(150,35){\circle*{2}}
\put(29.3,5.7){\vector(-1,1){28.6}}
\put(120.7,5.7){\vector(1,1){28.6}}
\put(120,6){\vector(0,1){28}}
\put(90.7,5.7){\vector(1,1){28.6}}
\put(90,26){\vector(0,1){6}}
\put(90,22){\line(0,-1){6}}
\put(90,13){\line(0,-1){6}}
\put(74.5,20.5){\vector(-1,1){13.8}}
\put(75.5,19.5){\line(1,-1){13.9}}
\put(60.7,5.7){\vector(1,1){28.6}}
\put(61,35){\vector(1,0){28}}
\put(59.5,14){\vector(0,-1){7}}
\put(60.5,14){\vector(0,-1){7}}
\put(59.5,17){\line(0,1){7}}
\put(60.5,17){\line(0,1){7}}
\put(59.5,27){\line(0,1){6}}
\put(60.5,27){\line(0,1){6}}
\put(59,35){\vector(-1,0){28}}
\put(59.3,5.7){\vector(-1,1){28.6}}
\bezier{400}(60,35)(20,45)(2,36)
\put(2,36){\vector(-3,-2){1}}
\put(89,5){\vector(-1,0){28}}
\put(59.2,5.4){\vector(-2,1){58.4}}
\end{picture}
\caption{$\Phi_3$, no maxima; 89320 (left) and 75108 (right)}
\label{75108}
\end{figure}

\begin{figure}
\unitlength 0.45mm
\begin{picture}(150,50)
\put(30,5){\circle{2}}
\put(60,5){\circle{2}}
\put(90,5){\circle{2}}
\put(120,5){\circle{2}}
\put(0,35){\circle*{2}}
\put(30,35){\circle*{2}}
\put(60,35){\circle*{2}}
\put(90,35){\circle*{2}}
\put(120,35){\circle*{2}}
\put(150,35){\circle*{2}}
\put(29.3,5.7){\vector(-1,1){28.6}}
\put(31,5){\vector(1,0){28}}
\put(59.3,5.7){\vector(-1,1){28.6}}
\put(59,35){\vector(-1,0){28}}
\put(59.5,14){\vector(0,-1){7}}
\put(60.5,14){\vector(0,-1){7}}
\put(59.5,17){\line(0,1){7}}
\put(60.5,17){\line(0,1){7}}
\put(59.5,27){\line(0,1){6}}
\put(60.5,27){\line(0,1){6}}
\put(45.5,20.5){\vector(1,1){13.8}}
\put(44.5,19.5){\line(-1,-1){13.8}}
\put(75.5,20.5){\vector(1,1){13.8}}
\put(74.5,19.5){\line(-1,-1){13.8}}
\put(89.3,5.7){\vector(-1,1){28.6}}
\put(61,35){\vector(1,0){28}}
\put(89,5){\vector(-1,0){28}}
\put(90,26){\vector(0,1){6}}
\put(90,22){\line(0,-1){6}}
\put(90,13){\line(0,-1){6}}
\put(120.7,5.7){\vector(1,1){28.6}}
\put(119.3,5.7){\vector(-1,1){28.6}}
\put(105.5,20.5){\line(1,1){13.8}}
\put(104.5,19.5){\line(-1,-1){13.8}}
\end{picture} \qquad 
\begin{picture}(150,50)
\put(30,5){\circle{2}}
\put(60,5){\circle{2}}
\put(90,5){\circle{2}}
\put(120,5){\circle{2}}
\put(0,35){\circle*{2}}
\put(30,35){\circle*{2}}
\put(60,35){\circle*{2}}
\put(90,35){\circle*{2}}
\put(120,35){\circle*{2}}
\put(150,35){\circle*{2}}
\put(59.5,14){\vector(0,-1){7}}
\put(60.5,14){\vector(0,-1){7}}
\put(59.5,17){\line(0,1){7}}
\put(60.5,17){\line(0,1){7}}
\put(59.5,27){\line(0,1){6}}
\put(60.5,27){\line(0,1){6}}
\put(31,5){\vector(1,0){28}}
\put(59.2,5.4){\vector(-2,1){58.4}}
\put(59.3,5.7){\vector(-1,1){28.6}}
\put(59,35){\vector(-1,0){28}}
\bezier{400}(2,36)(20,45)(60,35)
\put(2,36){\vector(-2,-1){1}}
\put(13,22){\vector(-1,1){9}}
\put(26,9){\line(-1,1){9}}
\put(61,35){\vector(1,0){28}}
\put(90,6){\vector(0,1){28}}
\put(120.7,5.7){\vector(1,1){28.6}}
\put(120,6){\vector(0,1){28}}
\put(90.7,5.7){\vector(1,1){28.6}}
\put(90,6){\vector(0,1){28}}
\put(60.7,5.7){\vector(1,1){28.6}}
\put(45.5,20.5){\vector(1,1){13.8}}
\bezier{100}(44.5,19.5)(42.5,17.5)(40.5,15.5)
\put(39.5,14.5){\line(-1,-1){8.5}}
\end{picture}
\caption{$\Phi_3$, no maxima; 63756 (left) and 59862 (right)}
\label{59862}
\end{figure}

\begin{figure}
\unitlength 0.45mm
\begin{picture}(150,50)
\put(30,5){\circle{2}}
\put(60,5){\circle{2}}
\put(90,5){\circle{2}}
\put(120,5){\circle{2}}
\put(0,35){\circle*{2}}
\put(30,35){\circle*{2}}
\put(60,35){\circle*{2}}
\put(90,35){\circle*{2}}
\put(120,35){\circle*{2}}
\put(150,35){\circle*{2}}
\put(89.5,14){\vector(0,-1){7}}
\put(90.5,14){\vector(0,-1){7}}
\put(89.5,17){\line(0,1){7}}
\put(90.5,17){\line(0,1){7}}
\put(89.5,27){\line(0,1){6}}
\put(90.5,27){\line(0,1){6}}
\put(29.3,5.7){\vector(-1,1){28.6}}
\put(30,6){\vector(0,1){28}}
\put(59.3,5.7){\vector(-1,1){28.6}}
\put(60,26){\vector(0,1){6}}
\put(60,22){\line(0,-1){6}}
\put(60,13){\line(0,-1){6}}
\put(89.3,5.7){\vector(-1,1){28.6}}
\put(89,35){\vector(-1,0){28}}
\put(75.5,20.5){\vector(1,1){13.8}}
\put(74.5,19.5){\line(-1,-1){13.8}}
\put(105.5,20.5){\vector(1,1){13.8}}
\put(104.5,19.5){\line(-1,-1){13.8}}
\put(119.3,5.7){\vector(-1,1){28.6}}
\put(119,5){\vector(-1,0){28}}
\put(120.7,5.7){\vector(1,1){28.6}}
\put(91,35){\vector(1,0){28}}
\put(61,5){\vector(1,0){28}}
\end{picture}
\qquad
\begin{picture}(150,50)
\put(30,5){\circle{2}}
\put(60,5){\circle{2}}
\put(90,5){\circle{2}}
\put(120,5){\circle{2}}
\put(0,35){\circle*{2}}
\put(30,35){\circle*{2}}
\put(60,35){\circle*{2}}
\put(90,35){\circle*{2}}
\put(120,35){\circle*{2}}
\put(150,35){\circle*{2}}
\put(89.5,14){\vector(0,-1){7}}
\put(90.5,14){\vector(0,-1){7}}
\put(89.5,17){\line(0,1){7}}
\put(90.5,17){\line(0,1){7}}
\put(89.5,27){\line(0,1){6}}
\put(90.5,27){\line(0,1){6}}
\put(60,26){\vector(0,1){6}}
\put(60,22){\line(0,-1){6}}
\put(60,13){\line(0,-1){6}}
\put(29.3,5.7){\vector(-1,1){28.6}}
\put(45.5,20.5){\vector(1,1){13.8}}
\put(44.5,19.5){\line(-1,-1){13.8}}
\put(59.3,5.7){\vector(-1,1){28.6}}
\put(61,5){\vector(1,0){28}}
\put(119,5){\vector(-1,0){28}}
\put(89,35){\vector(-1,0){28}}
\put(91,35){\vector(1,0){28}}
\put(119.3,5.7){\vector(-1,1){28.6}}
\put(105.5,20.5){\vector(1,1){13.8}}
\put(104.5,19.5){\line(-1,-1){13.8}}
\put(137,22){\vector(1,1){9}}
\put(123,8){\line(1,1){9}}
\put(75.5,20.5){\vector(1,1){13.8}}
\put(74.5,19.5){\line(-1,-1){13.8}}
\put(89.3,5.7){\vector(-1,1){28.6}}
\bezier{400}(148,36)(130,45)(90,35)
\put(148,36){\vector(2,-1){1}}
\put(110.5,15.25){\vector(2,1){38.5}}
\put(109.5,14.75){\line(-2,-1){18.5}}
\end{picture}
\caption{$\Phi_3$, no maxima; 53130 (left) and 34870 (right)}
\label{34870}
\end{figure}

\begin{figure}
\unitlength 0.45mm
\begin{picture}(150,50)
\put(30,5){\circle{2}}
\put(60,5){\circle{2}}
\put(90,5){\circle{2}}
\put(120,5){\circle{2}}
\put(0,35){\circle*{2}}
\put(30,35){\circle*{2}}
\put(60,35){\circle*{2}}
\put(90,35){\circle*{2}}
\put(120,35){\circle*{2}}
\put(150,35){\circle*{2}}
\put(89.5,14){\vector(0,-1){7}}
\put(90.5,14){\vector(0,-1){7}}
\put(89.5,17){\line(0,1){7}}
\put(90.5,17){\line(0,1){7}}
\put(89.5,27){\line(0,1){6}}
\put(90.5,27){\line(0,1){6}}
\put(60,26){\vector(0,1){6}}
\put(60,22){\line(0,-1){6}}
\put(60,13){\line(0,-1){6}}
\put(29.3,5.7){\vector(-1,1){28.6}}
\put(30,6){\vector(0,1){28}}
\put(59,5.7){\vector(-1,1){28.6}}
\put(61,5){\vector(1,0){28}}
\put(119,5){\vector(-1,0){28}}
\put(89,35){\vector(-1,0){28}}
\put(91,35){\vector(1,0){28}}
\put(119.3,5.7){\vector(-1,1){28}}
\put(105.5,20.5){\vector(1,1){13.8}}
\put(104.5,19.5){\line(-1,-1){13.8}}
\put(137,22){\vector(1,1){9}}
\put(123,8){\line(1,1){9}}
\put(75.5,20.5){\vector(1,1){13.8}}
\put(74.5,19.5){\line(-1,-1){13.8}}
\put(89.3,5.7){\vector(-1,1){28.6}}
\bezier{400}(90,35)(130,45)(148,36)
\put(148,36){\vector(2,-1){1}}
\put(110.5,15.25){\vector(2,1){38.4}}
\put(109.5,14.75){\line(-2,-1){18.5}}
\end{picture}
\caption{$\Phi_3$, no maxima; 29370}
\label{29370}
\end{figure}

\FloatBarrier

\begin{theorem}
\label{cher}
1. All virtual components of the type $\Phi_1$ with ten real critical points are chiral except for those with exactly two maxima or two minima $($see Fig.~\ref{26378}$)$, which are achiral. 
All virtual components of the type $\Phi_1$ with less than ten real critical points are achiral except for the ones with eight real critical points, having no minima or no maxima. 

2. All virtual components of the type $\Phi_3$ with ten real critical points are chiral except for those with exactly one maximum or one minimum $($see Fig.~\ref{77374}$)$, which are achiral. 
All virtual components of the type $\Phi_3$ with less than ten real critical points are achiral. 
\end{theorem}

\subsection{On the proofs of enumeration theorems  \ref{enu1} and \ref{enu3}}
\label{onproofs}
All possible values of the $\mbox{Card}$ and D-graph invariants mentioned in Theorems \ref{enu1} and \ref{enu3}
were found using the computer program described in \cite{AGLV2}, \cite{VS}.
For any type $\Phi_1$ or $\Phi_3$, we use the Gusein-Zade--A'Campo method (see \cite{AC}, \cite{GZ}) to compute the intersection matrix of vanishing cycles of a Morse polynomial of this class having only real critical points. Then, using theorem 1.4 of \S 5.1 in \cite{APLT}, we calculate the intersection indices of these vanishing cycles with the set of real points. Thus, we obtain a virtual Morse function associated with this polynomial.
Starting from this initial data, the combinatorial program \label{progg}

{\footnotesize\begin{verbatim}https://drive.google.com/file/d/19FK3NDqHr01CVyVegO6wNUaedrvi_8nJ/view?usp=sharing
\end{verbatim}}

\noindent
runs through the entire formal graph of this type and counts, in particular, 
the number of all virtual Morse functions of this type and the number of virtual functions with each passport invariant value. For each passport value for which this number is not zero, it then (upon request) finds a virtual Morse function with this value. Starting from this virtual Morse function, a slightly modified version 
of the same program (with virtual surgeries of types $s1$ and $s3$ disabled)
counts the number of virtual Morse functions in its virtual component, i.e., the Card invariant of this virtual function. (For example, this program 
{\footnotesize \begin{verbatim}https://drive.google.com/file/d/1bGV16NMqm-VxvChKxjLTqxORUtKeIVQe/view?usp=sharing
\end{verbatim}} 

\noindent 
with certain initial data handles the virtual component of the $\Phi_3$ class, which has eight real critical points and Card = 131148.)
If this number is smaller than the total number of virtual functions with that passport, we find another virtual function with the same passport but in a different virtual component. Then we calculate its $\mbox{Card}$ value as well. We continue this process until the sum of the different values of the Card invariant of virtual functions with any given passport reaches the total number of virtual functions with that passport. This computation (and the subsequent reconstruction of D-graphs from the virtual Morse functions representing these components in the case of polynomials with ten real critical points) proves all statements of Theorems \ref{enu1} and \ref{enu3}. 

In sections \ref{reali} and \ref{chirali}, we investigate the chirality of all these components and prove Theorem \ref{cher}. Specifically, in \S~\ref{reali} we realize nearly all achiral virtual components by polynomials that are symmetric with respect to the reflection (\ref{invol0}). In \S~\ref{achirali}, we describe a homological criterion of chirality in the terms of the formal graph and use it to demonstrate the achirality of the remaining achiral components. In \S~\ref{cchirali}, we prove the chirality of all chiral components.

\begin{remark} \rm The spaces $\Phi_1$ and $\Phi_3$ are invariant under the involution (\ref{invol}). This involution maps Morse polynomials with any value $(M, m_+)$ of the ``passport'' invariant to polynomials with the value $(M, M/2-m_+)$ in the case $\Phi_1$ and to polynomials with the value $(M, M/2 -1 -m_+) $ in the case $\Phi_3$. Therefore, it is sufficient to study only the virtual Morse functions with $m_+ \leq M/4$ (respectively, $m_+ \leq M/4 -1$).
\end{remark}

\subsection{Normalization of $D$-graphs}
\label{norma}

\begin{definition} \rm
An edge of a $D$-graph is called {\em normal} if 

a) it is oriented from a vertex with a smaller Morse index to a vertex with a larger index,

b) it is dashed if the parities of Morse indices of the critical points corresponding to its ends are the same; it is solid if these parities are different. 

 Otherwise, this edge is called a {\em tunnel} edge. The {\em normalization} of a $D$-graph consists in removing all its tunnel edges.
\end{definition}

By analyzing Figs.~\ref{122298}--\ref{29370}, we see that

a) all $D$-graphs of polynomials with both maxima and minima are already normal,

b) seventeen out of eighteen $D$-graphs of polynomials without maxima shown in these figures are split by the normalization into pairs of standard Coxeter--Dynkin graphs of some simple singularities after normalization;

c) the remaining $D$-graph (see Fig.~\ref{102234} left) splits into an isolated vertex and the extended Coxeter--Dynkin graph of type $\tilde E_8$.
\medskip

A very similar situation occurs for $D$-graphs of other parabolic singularities, see \cite{Vx9} and \cite{Vp8} (with extended graphs of types $\tilde E_7$ and $\tilde E_6$).

\section{Proof of Proposition \ref{uniq}}
\label{secproof}

\begin{lemma}
\label{lempre}
The subset of the space $\Phi_1$, which consists of polynomials having a critical point of class $E_8$ with zero critical value, and also a critical point of class $ A_2$, is the union of four smooth components diffeomorphic to $\R^7$ 
and moved to each other by the involutions
 $($\ref{invol0}$)$ and $($\ref{invol}$)$. Two of these components consist of polynomials with negative critical values at the $A_2$ points, and the other two consist of polynomials with positive values. Each component is swept out by a one-parametric family of orbits of the group $\G$. Each of these orbits contains exactly one polynomial of the form 
\begin{equation}
x^3 \pm y^6 + \tau y^5 \label{normaz} , \quad \tau \neq 0 \ . 
\end{equation} 
The intersection of each such component with the parameter space of the deformation $($\ref{vers1}$)$ is an orbit of the group $\R_+$ of positive quasihomogeneous dilations \begin{equation}
\label{dila} T_t: f(x, y) \to t^{-6} f(t^2x, t y), \quad t > 0. \end{equation}
\end{lemma}

\noindent
{\it Proof.} Let $f $ be a polynomial of type $\Phi_1$ having a real critical point of class $E_8$ with zero critical value. Its 2-jet at this point is trivial, and 3-jet has the form $q(x-\nu y)^3,$ $q>0$. 
The substitution $\tilde x = x - \nu y$ and a dilation of $x$ coordinate reduce this 3-jet to the form $\tilde x^3$. The resulting polynomial $\tilde f \in \Phi_1$ has zero coefficients at the monomials $y^4$ and $\tilde x y^3$, because otherwise its critical point at the origin would be of class $E_6$ or $E_7$. 
Conversely, the coefficient at the monomial $y^5$ is non-zero: otherwise the Milnor number is greater than $8$. So, this polynomial is a linear combination of monomials $\tilde x^3,$ $\tilde x^2y^2$, $\tilde x y^4$, $y^6$, and $y^5$.
An additional substitution of the form $\check x = \tilde x - \eta y^2$ and dilation of $y$ coordinate (which belong to the group $\G$) reduce $\tilde f$ to the form
\begin{equation} 
\label{alter} \check x^3 + a \check x \check y^4 \pm \check y^6 + \tau \check y^5 \qquad \mbox{or} \qquad \check x^3 \pm \check x \check y^4 + \tau \check y^5, 
\end{equation} 
$\tau \neq 0.$ 
The condition of having also a critical point of class $A_2$ (that is, the coincidence of two non-zero complex solutions of the system $f'_x =0 =f'_y$) prohibits the right-hand possibility of the alternative (\ref{alter}) and implies the condition $a=0$ in the left-hand possibility. The resulting set of polynomials of the form (\ref{normaz}) 
consists of four components, which are characterized by the sign before the monomial $y^6$ and the sign of the coefficient $\tau$. 

All points (\ref{normaz}) belong to different $\G$-orbits. Indeed, the critical value at the $A_2$-point is an invariant of orbits. This value is a monomial of degree six of the parameter $\tau$ in (\ref{normaz}). In particular, for any $c \neq 0$ there are exactly two polynomials (\ref{normaz}) having the critical value equal to $c$. These two 
two polynomials are obtained from each other by the multiplication of $\tau$ by $-1$. They are not $\G$-equivalent because the $y$ coordinates of their $A_2$ points have different signs. So, the subset of $\Phi_1$ considered in Lemma is the family of $\G$-orbits represented and parametrized by polynomials (\ref{normaz}).
 The group $\Z_2^2$ of two involutions (\ref{invol0}), (\ref{invol}) maps all its four components into each other. If a polynomial of the form (\ref{vers1}) belongs to one of these components, then its entire orbit under the action of the one-parametric group 
(\ref{dila}) also belongs to it. These one-dimensional orbits are smooth curves in the space of polynomials (\ref{vers1}), and the $\G$-orbits intersect this space transversally. Therefore, the four components in consideration are canonically diffeomorphic to the products of these smooth curves and the group $\G$. \hfill $\Box$

\begin{lemma}
\label{lempost}
The subset of the space $\Phi_3$, which consists of polynomials having a critical point of class $D_8^-$ with zero critical value, and also a class $ A_2$ critical point, is the union of four connected components that are diffeomorphic to $\R^7$ 
and are moved to each other by the involutions 
 $($\ref{invol0}$)$ and $($\ref{invol}$)$.
Two of these components consist of polynomials with the negative critical values at the $A_2$ critical points, and the other two consist of polynomials with the positive critical values. Each of these four components is swept out by a one-parametric family of orbits of the group $\G$. Each of these orbits contains exactly one polynomial of the form 
\begin{equation}
\label{normaw}
x^3\pm \sqrt{21}x^2y^2 -x y^4 + \alpha x^2y \ .
\end{equation}
The intersection of each component with the parameter space of deformation $($\ref{vers3}$)$ is an orbit of the group $\R_+$ of positive quasihomogeneous dilations $($\ref{dila}$)$.
\end{lemma}

\noindent
{\it Proof.}
Each polynomial $f$ of type $\Phi_3$ with a critical point of class $D_8^-$ and critical value 0 at this point can be reduced by the action of the group $\G$ to the form 
\begin{equation}
\label{ee}
x^3 + b x^2 y^2 \pm x y^4 + \alpha x^2 y, \quad \alpha \neq 0 . 
\end{equation}
Indeed, the 2-jet of $f$ at its $D_8^-$ critical point is trivial, and the 3-jet has the form $q (x - \varkappa y)^2(x - \nu y)$ for some $q >0$ and $\varkappa \neq \nu$. The substitution $\tilde x = x - \varkappa y$ and a dilation of $\tilde x$ coordinate reduce this 3-jet to the form $\tilde x^3 + \alpha \tilde x^2 y$, $\alpha \neq 0$.
The resulting polynomial has zero coefficient at the monomial $y^4$: otherwise this polynomial would be of class $D_5$. 
Its lower quasihomogeneous part with weights $\deg \tilde x = 2, \deg y=1$
should have the form $y (\tilde x - \theta y^2)^2$, because otherwise the polynomial is of the class $D_6$. So, this quasihomogeneous part can be reduced to the monomial $\tilde x^2 y$ by a diffeomorphism of class $\G$ keeping the entire polynomial within the space $\Phi_3$. The coefficient of the resulting polynomial at the monomial $y^6$ should vanish, because otherwise the critical point is of class $D_7$. Finally, the coefficient at the monomial $\tilde x y^4$ is not zero, because otherwise the Milnor number would be greater than $8$. Dilating the coordinate $y$, we get a polynomial of the form (\ref{ee}).

The sum of Milnor numbers of all critical points in $\C^2$ of the obtained polynomial is equal to ten, so there are only two critical points outside the origin (counting with the multiplicities). The system of equations $f'_x =0 =f'_y$ for these two points implies easily the condition 
$$ \pm 48 y^2 = 12 b^2 y^2 + 20 b \alpha y + 7 \alpha^2$$
on their $y$-coordinate, where the sign at $48 y^2$ repeats the sign at $x y^4$ in (\ref{ee}).
If these two critical points coincide (at a point of class $A_2$), then $b^2 \pm 21=0$, so the sign $\pm $ in (\ref{ee}) is equal to $-$ and $b $ is equal to $\sqrt{21}$ or $-\sqrt{21}$. The critical values of the resulting polynomials 
(\ref{normaw}) at the $A_2$ point
 depend on the parameter $\alpha$ as a monomial of degree six. They are negative in the case of the sign $+$ at $\sqrt{21}x^2y^2$ and are positive in the case of the sign $-$.
The concluding arguments are the same as in the proof of Lemma \ref{lempre}. 
 \hfill $\Box$

\begin{corollary}
\label{lemcor}
1. The space of polynomials $($\ref{vers1}$)$ contains exactly two polynomials having a critical point of class $E_8$ with zero critical value and a critical point of class $A_2$ with critical value $-\frac{1}{5}$. Namely, these two polynomials are the intersection points of this space with the $\G$-orbits of the polynomials $x^3+y^6 \pm \frac{6}{5}y^5$.

2. The space of polynomials $($\ref{vers3}$)$ contains exactly two polynomials having a critical point of class $D_8^-$ with zero critical value and a critical point of class $A_2$ with critical value $-\frac{343343}{375000}\sqrt{\frac{7}{3}}$. Namely, these two polynomials are the intersection points of this space and the $\G$-orbits of polynomials $x^3 + \sqrt{21} x^2y^2-x y^4 \pm 6 x^2y.$

In both cases, the two polynomials are mapped to each other by the involution $($\ref{invol0}$)$. \hfill $\Box$
\end{corollary}

\noindent
{\it Proof.} Proof follows immediately from Lemmas \ref{lempre} and \ref{lempost}, Proposition \ref{proponne}, and direct calculations with polynomials (\ref{normaz}) and (\ref{normaw}).
\medskip

\begin{figure}
\unitlength 0.8mm
\begin{picture}(70,36)
\bezier{70}(0,23)(5,20)(7,20)
\bezier{70}(7,20)(11,19)(20.2,19)
\bezier{100}(20.2,19)(30,20)(40,17)
\bezier{100}(40,17)(50,10)(60,11)
\bezier{70}(60,11)(70,14)(60,18)
\bezier{100}(60,18)(55,20)(45,17)
\bezier{150}(45,17)(25,10)(5,31)
\bezier{100}(5,31)(1,39)(7,36)
\bezier{100}(7,36)(10,33)(13,20)
\bezier{100}(13,20)(14,17)(20,14)
\bezier{150}(20,14)(40,5)(60,0)
\put(12,24.5){\circle*{1}}
\put(13.5,19){\circle*{1}}
\put(20.2,19){\circle*{1}}
\put(41.5,16){\circle*{1}}
\put(13,4){\footnotesize $-$}
\put(5.5,30){\footnotesize $-$}
\put(13.5,19.2){\scriptsize $-$}
\put(29,15.6){\footnotesize $-$}
\put(52,14){\footnotesize $-$}
\put(33,27){\footnotesize $+$}
\end{picture} \qquad \quad
\begin{picture}(70,36)
\bezier{200}(0,14)(12,30)(28,18)
\bezier{200}(28,18)(45,6)(70,30)
\bezier{300}(0,22)(12,6)(28,18)
\bezier{300}(28,18)(45,30)(70,6)
\put(64,0){\line(0,1){36}}
\put(14,16.5){\scriptsize $-$}
\put(39,16.5){\scriptsize $-$}
\put(58.5,16.5){\scriptsize $-$}
\put(3.7,18){\circle*{1}}
\put(28,18){\circle*{1}}
\put(54.2,18){\circle*{1}}
\put(64,11.33){\circle*{1}}
\put(64,24.67){\circle*{1}}
\end{picture}
\caption{Morse perturbations of $E_8$ and $D_8^-$ singularities}
\label{J38}
\end{figure}

\begin{lemma}
\label{le39}
Each of the two polynomials considered in statement 1 of Corollary \ref{lemcor} has arbitrarily small Morse perturbations which 

a$)$ split its critical point of class $E_8$ into four saddlepoints with zero critical value and four minima in such a way that the local zero-level set of the obtained polynomial
looks as shown in Fig.~\ref{J38} $($left$)$ or its reflection about a vertical line, and 

b$)$ split the $A_2$ critical point into one minimum point and one saddlepoint with critical values near $-\frac{1}{5}$.

The $D$-graph of these Morse perturbations is shown in Fig.~\ref{3852} $($left$)$.
\end{lemma}

\begin{lemma}
\label{le39a}
Each of the two polynomials considered in statement 2 of Corollary \ref{lemcor}
 has arbitrarily small Morse perturbations which 

a$)$ split its critical point of class $D_8^-$ into five saddlepoints with zero critical value and three minima in such a way that the local zero-level set of the obtained polynomial looks as shown in Fig.~\ref{J38} $($right$)$ or its reflection with respect to a vertical line, and 

b$)$ split the $A_2$ critical point into one minimum point and one saddlepoint with critical values near \ $-\frac{343343}{375000}\sqrt{\frac{7}{3}}$.

The $D$-graph of these Morse perturbations is shown in Fig.~\ref{59862} $($right$)$.
\end{lemma}

\noindent
{\it Proofs of Lemmas \ref{le39}, \ref{le39a}.} 
The existence of such independent perturbations of two critical points follows from the versality of the deformations (\ref{vers1}), (\ref{vers3}) and from the perturbations of $E_8$ and $D_8^-$ singularities demonstrated in pp. 17 and 15 of \cite{AC}.
By the construction, the vertices of the D-graph of the first (respectively, the second) perturbation can be split into subsets of cardinality eight and two in such a way that the corresponding subgraphs are the canonical Coxeter-Dynkin graphs of types $E_8$ (respectively, $D_8$) and $A_2$. Among all D-graphs listed in Theorem \ref{enu1} (respectively, \ref{enu3}), only the D-graph called in the last assertion of Lemma \ref{le39} (respectively, Lemma \ref{le39a}) allows such a splitting. \hfill $\Box$ \medskip

The proof of Proposition \ref{uniq} follows from the statements \ref{lempre}--\ref{le39a} in the same way as Proposition 3 of \cite{Vx9} is deduced from Lemmas 4 and 5 and Proposition 22. First, we describe this proof for the type $\Phi_1$ polynomials, referring to \cite{Vx9} for the common details.

Using the Lyashko--Looijenga covering, we can assume that our two Morse polynomials $f$ and $ \tilde f$ have equal sets of critical values and equal systems of paths in $\C^1$ defining the vanishing cycles. Using the diffeomorphisms of class $\G$ we can also assume that they lie in the space of deformation (\ref{vers1}). We draw a generic piecewise-algebraic path $I: [0,1] \to \Phi_1$ in this space, connecting the polynomial $f \equiv I(0)$ with one of two points of the $(E_8, A_2)$ stratum considered in Corollary \ref{lemcor}. Using the Lyashko--Looijenga map and the coincidence of virtual Morse functions associated with $f$ and $\tilde f$, we can uniquely draw a path $\tilde I:[0,1] \to \Phi_1$ starting from the polynomial $\tilde f$, repeating the sets of critical values of the functions of the first path and undergoing the same standard surgeries. The endpoint $\tilde I(1)$ of this path has the critical points of classes $E_8$ and $A_2$ with the same critical values as $I(1)$, so according to Corollary \ref{lemcor}
it either coincides with the point $I(1)$ or is symmetric to it via the involution (\ref{invol0}). 

In the second case, when $I(1) \neq \tilde I(1)$, we apply this involution to entire path $\tilde I$, including its starting point $\tilde f$, and reduce the proof to the first case.

In the first case, the final segments of the paths $I$ and $\tilde I$ should coincide in a neighborhood of their endpoint $I(1) = \tilde I(1)$. Indeed, suppose that the points $I(1-\varepsilon)$ and $\tilde I(1-\varepsilon)$ are different for small $\varepsilon>0$. We can assume that these points lie in the set of generic polynomials. Consider the map from one of these points to the other. Since their associated virtual Morse functions and critical values are the same, the covering homotopy of the Lyashko--Looijenga map extends this maps to an automorphism of the space $\Phi_1$, sending the point $I(1)$ to itself and commuting with the Lyashko--Looijenga map. 

Analyzing the multisingularity $(E_8, A_2)$ as in \S~3.5 of \cite{Vx9} proves that such an automorphism can only be the identity. Indeed, by the versality property, a neighborhood of the point $I(1)$ is naturally diffeomorphic to the Cartesian product of the parameter spaces of the canonical versal deformations of the simple singularities of classes $E_8$ and $A_2$. Namely, in this neighborhood the family (\ref{vers1}) can be considered as a deformation of the $E_8$ singularity of the polynomial $I(1)$, and also as a deformation of $A_2$ singularity of the polynomial $I(1)+\frac{1}{5}$. The projection maps of the Cartesian product structure to its factors are maps of the parameter spaces that locally induce the deformation (\ref{vers1}) from these two deformations as in the definition of versal deformations, see \S~8 in  \cite{AVGZ}. The aforementioned automorphism of a neighborhood of the point $I(1)$ then induces the automorphisms of the canonical versal deformations of the $E_8$ and $A_2$ singularities that commute with the corresponding Lyashko--Looijenga coverings. All complex automorphisms of versal deformations of simple singularities that commute with the Lyashko--Looijenga covering are listed in \cite{Liv}. This list implies that there are no such non-trivial automorphisms that preserve the space of real functions. Thus, the paths $I$ and $\tilde I$ coincide near their common endpoint, and hence they coincide everywhere, including their starting points $f$ and $\tilde f$. This proves Proposition \ref{uniq} for polynomials of class $\Phi_1$.

In the case $\Phi_3$, an additional difficulty occurs at the final part of the proof. Namely,
the canonical real versal deformation 
\begin{equation}
\label{vd8}
 F(\xi, \lambda) \equiv \xi^2\eta -\eta^7 +\lambda_1 +\lambda_2 \xi + \lambda_3 \eta + \lambda_4 \eta^2 + \lambda_5 \eta^3 + \lambda_6 \eta^4 + \lambda_7 \eta^5 + \lambda_8 \eta^6
\end{equation}
of $D_8^-$ singularity has exactly one non-trivial automorphism commuting with the Lyashko--Looijenga map: it is defined by the reflection 
\begin{equation}
\label{refX}
(\xi, \eta) \leftrightarrow (-\xi, \eta).
\end{equation}
This fact follows from the description of all complex automorphisms of the $D_8$ singularities that commute with the Lyashko--Looijenga covering (see \cite{Liv}), and from the selection of automorphisms that preserve the set of real polynomials.

Therefore, a neighborhood of the point $I(1)$ in the space of polynomials (\ref{vers3}) has exactly one such automorphism as well.
Specifically, this neighborhood can be naturally identified with the direct product of the parameter spaces of the standard versal deformations of the $D_8^-$ and $A_2$ singularities. The unique nontrivial automorphism of the neighborhood of the point $I(1)$ acts as the product of the automorphism  on the eight-dimensional factor induced from by the reflection (\ref{refX}) and the identity automorphism of the versal deformation of the $A_2$ singularity.

Let us prove that this unique nontrivial local automorphism does not preserve the associated virtual Morse functions.

Denote by $U$ the set of polynomials of the form (\ref{vers3}) near the point $I(1)$, such that the $D_8^-$-point of the polynomial $I(1)$ splits into three minima with negative critical values and five saddlepoints with positive critical values, and the $A_2$-point splits into two real critical points, a minimum and a saddlepoints with critical values near $-\frac{1}{5}$. According to Theorem 1.7 of \cite{Lo0} and the direct product structure described above, this set consists of a single contractible connected component. We can assume that the path $I$ approaches the point $I(1)$ from inside this domain $U$.

\begin{figure}
\unitlength 0.8mm
\begin{picture}(70,36)
\bezier{100}(0,15.5)(3,18)(0,20.5)
\bezier{200}(9,14.5)(1,18)(9,21.5)
\bezier{200}(22,14.5)(30,18)(22,21.5)
\bezier{200}(22,14.5)(15.5,12)(9,14.5)
\bezier{200}(22,21.5)(15.5,24)(9,21.5)
\bezier{200}(34,14.5)(26,18)(34,21.5)
\bezier{200}(47,14.5)(55,18)(47,21.5)
\bezier{200}(47,14.5)(40.5,12)(34,14.5)
\bezier{200}(47,21.5)(40.5,24)(34,21.5)
\bezier{100}(60,14.5)(53,18)(60,21.5)
\bezier{100}(60,21.5)(64,25)(64,18)
\bezier{100}(60,14.5)(64,13)(64,18)
\bezier{150}(64,0)(64,11.5)(70,6)
\bezier{150}(64,36)(64,22)(70,30)
\put(14,16.5){\scriptsize $-$}
\put(39,16.5){\scriptsize $-$}
\put(58.5,16.5){\scriptsize $-$}
\put(3.3,18){\circle*{1}}
\put(28,18){\circle*{1}}
\put(54.2,18){\circle*{1}}
\put(64,11.33){\circle*{1}}
\put(64,24.67){\circle*{1}}
\put(3,14){\tiny 6}
\put(27,14){\tiny 7}
\put(53,14){\tiny 8}
\put(65,22.5){\tiny 9}
\put(65,12){\tiny 10}
\end{picture}
\qquad \qquad
\begin{picture}(70,36)
\bezier{100}(0,15.5)(3,18)(0,20.5)
\bezier{200}(9,14.5)(1,18)(9,21.5)
\bezier{200}(22,14.5)(30,18)(22,21.5)
\bezier{200}(22,14.5)(15.5,12)(9,14.5)
\bezier{200}(22,21.5)(15.5,24)(9,21.5)
\bezier{200}(34,14.5)(26,18)(34,21.5)
\bezier{200}(47,14.5)(55,18)(47,21.5)
\bezier{200}(47,14.5)(40.5,12)(34,14.5)
\bezier{200}(47,21.5)(40.5,24)(34,21.5)
\bezier{100}(60,14.5)(53,18)(60,21.5)
\bezier{100}(60,21.5)(64,23)(64,18)
\bezier{100}(60,14.5)(64,11)(64,18)
\bezier{150}(64,0)(64,14)(70,6)
\bezier{150}(64,36)(64,26)(70,30)
\put(14,16.5){\scriptsize $-$}
\put(39,16.5){\scriptsize $-$}
\put(58.5,16.5){\scriptsize $-$}
\put(3.3,18){\circle*{1}}
\put(28,18){\circle*{1}}
\put(54.2,18){\circle*{1}}
\put(64,11.33){\circle*{1}}
\put(64,24.67){\circle*{1}}
\put(3,14){\tiny 6}
\put(27,14){\tiny 7}
\put(53,14){\tiny 8}
\put(65,22.5){\tiny 10}
\put(65,12){\tiny 9}
\end{picture}
\caption{Ordering the critical values of perturbations of $D_8^-$ singularity}
\label{next}
\end{figure}

 The boundary of the domain $U$ contains
a Morse polynomial considered in Lemma \ref{le39a}. Therefore, the D-graphs of
all polynomials from $U$ also are as shown in Fig.~\ref{59862} (right).

The perturbations of $I(1)$ that satisfy the conditions of Lemma \ref{le39a} are not generic, as their five Morse points have common critical value 0. Such a perturbation $f_\lambda$ can be slightly perturbed again so that the resulting polynomial $f_{\lambda'}$ is generic and belongs to the set $U$, and the critical values at all its five saddlepoints 
obtained by splitting the $D_8^-$ singularity are ordered in $\R^1$ as shown by the numbers in Fig.~\ref{next} (left) at the corresponding saddlepoints of the polynomial $f_\lambda$. (These numbers start with 6 since the lower critical values at the other five critical points are taken into account.) Let $f_{\lambda''}$ be the Morse polynomial obtained from $f_{\lambda'}$ by our local automorphism. It is very close to $f_\lambda$ and hence to $f_{\lambda'}$.
The critical values of $f_{\lambda''}$ are ordered as shown in Fig.~\ref{next} (right).
According to the Gusein-Zade--A'Campo calculation method of intersection indices,
the vanishing cycles in the sets $f^{-1}_{\lambda'}(0)$ and $f^{-1}_{\lambda''}(0)$
that vanish at the critical points labeled by 9 and 10 correspond to the horns of the Coxeter--Dynkin graph \ \
 \begin{picture}(47,7)
\unitlength 0.7 mm
\put(0,4){\circle*{1.3}}
\put(20,4){\circle*{1.3}}
\put(40,4){\circle*{1.3}}
\put(58,0){\circle*{1.3}}
\put(58,8){\circle*{1.3}}
\put(10,4){\circle{1.3}}
\put(30,4){\circle{1.3}}
\put(50,4){\circle{1.3}}
\put(0,4){\line(1,0){9.5}}
\put(10.5,4){\line(1,0){19}}
\put(30.5,4){\line(1,0){19}}
\put(58,0){\line(-2,1){7.5}}
\put(58,8){\line(-2,-1){7.5}}
\put(-0.5,0){\tiny 6}
\put(19,0){\tiny 7}
\put(39,0){\tiny 8}
\put(59,-0.5){\tiny 10}
\put(59.4,7){\tiny 9}
\end{picture}
of the $D_8^-$ singularity.

Two generic perturbations $f_{\lambda'}$ and $f_{\lambda''}$ can be continuously deformed  into each other within the domain $U$ in such a way that the cycle vanishing at the point marked by 9 in Fig.~\ref{next} (left) will be moved to the cycle vanishing at the point marked by 10 in Fig.~\ref{next} (right), and the orders of critical values of all other critical points of $f_\lambda$ will not change. Therefore, if the perturbations $f_{\lambda'}$ and $f_{\lambda''}$ have equal associated virtual Morse functions, then the intersection indices of both vanishing cycles $\Delta_9$  and $\Delta_{10}$ in $f^{-1}_{\lambda'}(0)$ with all other basic vanishing cycles should be the same. However, Fig.~\ref{59862} (right) shows that the cycle corresponding to one horn of the $D_8$-subgraph of the D-graph has nonzero intersection indices with both cycles corresponding to the vertices of the $A_2$-subgraph, while the cycle corresponding to the other horn has  nonzero intersection index with only one of them. Thus, the virtual Morse functions associated with polynomials $f_{\lambda'}$ and $f_{\lambda''}$ are different.
Therefore, the paths $I$ and $\tilde I$ coincide near their endpoints and, consequently, everywhere, including their initial points $f$ and $\tilde f$. \hfill $\Box$

\section{Realization of virtual Morse functions}
\label{reali}

In this section, we realize by polynomials many virtual components, whose existence was predicted in Theorems \ref{enu1} and \ref{enu3}.
Almost all these realizations are invariant under the involution (\ref{invol0}) and hence prove the achirality of the corresponding isotopy classes. 
 
In all our pictures, it is assumed that the $y$ axis is horizontal and the $x$ axis is vertical and directed upwards.

\subsection{Realizations of virtual Morse functions of class $\Phi_1$}

A virtual Morse function with ten real critical points, of which exactly one is a local maximum, can be realized by the polynomial
\begin{equation}
(x^2 + y^4 - 8 y^2)(x+5(y+3/2)^2-5) .
\label{J10a}
\end{equation}
A virtual Morse function with ten real critical points, of which exactly two are local maxima, can be realized by the polynomial
\begin{equation}
(x^2 + y^4 - 2 y^2)(x+y^2-1), \label{J10b}
\end{equation}
in particular the corresponding connected component of the space of Morse polynomials is achiral.
The zero-level sets of these two polynomials are outlined in the left parts of Figs.~\ref{122298} and \ref{26378}. 

\begin{remark} \rm
The polynomials (\ref{J10a}) and (\ref{J10b}) are not generic since they have multiple critical value 0. However, their arbitrary small generic perturbations realize the promised virtual components, cf. Definition \ref{def8}.
\end{remark}

Virtual Morse functions with eight real critical points and exactly two (respectively, one, respectively, no) maximum points can be realized by polynomials whose zero-level sets are shown in Fig.~\ref{J18} left (respectively, Fig.~\ref{J18} right, respectively, Fig.~\ref{J38} left).

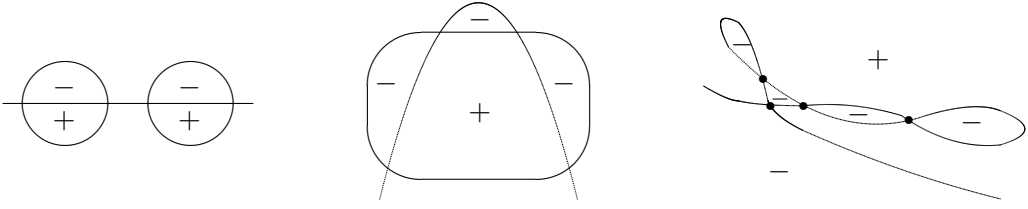
\begin{figure}
\unitlength 0.55mm
\begin{picture}(80,40)
\put(0,23){\line(1,0){60}}
\put(15,23){\circle{20}}
\put(45,23){\circle{20}}
\put(12,17){\footnotesize $+$}
\put(12,25){\footnotesize $-$}
\put(42,17){\footnotesize $+$}
\put(42,25){\footnotesize $-$}
\end{picture} \ \qquad \qquad 
\unitlength 0.55 mm
\begin{picture}(80,40)
\bezier{400}(15,0)(35,80)(55,0)
\put(35,19){\oval(45,30)}
\put(14,22){\footnotesize $-$}
\put(33,35){\footnotesize $-$}
\put(50,22){\footnotesize $-$}
\put(33,16){\footnotesize $+$}
\end{picture} 
\caption{Perturbations for $\Phi_1$ with eight real critical points}
\label{J18}
\end{figure}

For the first of them, we can take the polynomial 
\begin{equation}
x \left(x^2 + (y^2-1)^2 - \varepsilon \right), \quad \varepsilon \in (0, 1), 
\label{J8a}
\end{equation}
 For the second, we take the polynomial 
\begin{equation}
(x^2 + y^4 -1)(x+ 2y^2- A), \quad A \in ( 1 , \sqrt{2}). \label{J8b}
\end{equation}
To construct the third polynomial, we take the perturbation 
\begin{equation}
\label{e77}
x^3 + y^6 +\varepsilon (x^2 y+3y^5), \quad \varepsilon >0,
\end{equation} 
of singularity $x^3 + y^6 $. It has one critical point of class $E_8$ and no other real critical points. Then we apply the standard perturbation of this $E_8$ singularity, as described on page 17 of \cite{AC}. This perturbation can be chosen arbitrarily small, so that it does not return the non-real critical points of (\ref{e77}) to the real domain.

\begin{figure}
\unitlength 0.75mm
\begin{picture}(50,30)
\put(25.00,16.80){\circle*{1.33}}
\put(31.10,10.20){\circle*{1.33}}
\put(18.90,10.20){\circle*{1.33}}
\bezier{100}(1,1)(10,2)(17,8.5)
\bezier{300}(17,8.5)(40,34)(47,26)
\bezier{100}(47,26)(50,24)(47,20)
\bezier{400}(47,20)(25,-1)(3,20)
\bezier{100}(3,20)(0,24)(3,26)
\bezier{400}(3,26)(10,34)(33,8.5)
\bezier{100}(49,1)(40,2)(33,8.5)
\put(25.00,12.00){\small \makebox(0,0)[cc]{$+$}}
\put(39.00,20.00){\small \makebox(0,0)[cc]{$-$}}
\put(11.00,20.00){\small \makebox(0,0)[cc]{$-$}}
\put(25,0){\small \makebox(0,0)[cc]{$-$}}
\put(25,25){\small \makebox(0,0)[cc]{$+$}}
\end{picture} \qquad 
\begin{picture}(50,26)
\bezier{100}(1,1)(17,3)(23,10)
\bezier{300}(23,10)(40,31)(47,26)
\bezier{100}(47,26)(50,24)(47,22)
\bezier{320}(47,22)(25,3.5)(3,22)
\bezier{100}(3,22)(0,24)(3,26)
\bezier{300}(27,10)(10,31)(3,26)
\bezier{100}(27,10)(33,3)(49,1)
\put(25.00,12.50){\circle*{1.33}}
\put(40.00,22.00){\small \makebox(0,0)[cc]{$-$}}
\put(10.00,22.00){\small \makebox(0,0)[cc]{$-$}}
\put(25.00,0){\small \makebox(0,0)[cc]{$-$}}
\put(25.00,22){\small \makebox(0,0)[cc]{$+$}}
\end{picture} \qquad
\begin{picture}(50,26)
\bezier{100}(1,1)(17,3)(23,10)
\bezier{300}(23,10)(40,31)(47,26)
\bezier{100}(47,26)(50,24)(47,22)
\bezier{320}(47,22)(25,14)(3,22)
\bezier{100}(3,22)(0,24)(3,26)
\bezier{300}(27,10)(10,31)(3,26)
\bezier{100}(27,10)(33,3)(49,1)
\put(25.00,12.30){\circle*{1.33}}
\put(30,18.20){\circle*{1.33}}
\put(20,18.20){\circle*{1.33}}
\put(25.00,16.00){\small \makebox(0,0)[cc]{$-$}}
\put(40.00,22.00){\small \makebox(0,0)[cc]{$-$}}
\put(10.00,22.00){\small \makebox(0,0)[cc]{$-$}}
\put(22.5,2){\footnotesize $-$}
\put(22.5,24){\footnotesize $+$}
\end{picture}
\caption{Perturbations for $\Phi_1$ with six real critical points}
\label{J16}
\end{figure}
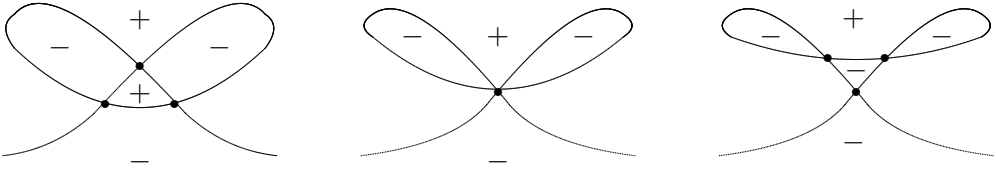

\medskip

The perturbation 
\begin{equation}
x^3 + y^6 - \varepsilon x y^2 , \quad \varepsilon >0, \label{J6a}
\end{equation}
 of $J_{10}^1 $
 singularity 
$x^3 + y^6$ has a critical point of class $D_4^-$, two local minima and no other critical points, see Fig.~\ref{J16} (center). Additional perturbations of its $D_4^-$ critical point, shown in Fig.~\ref{J16} (left and right), can be performed as indicated in \cite{AC} or \cite{Vsing} preserving the symmetry with respect to the coordinate $y$. These perturbations realize virtual Morse functions with six real critical points, exactly one (respectively, none) of which are maxima. 
The perturbation 
\begin{equation}
\label{J4b} x^3+ y^6 - \varepsilon x y^2 + \varepsilon^3 x
\end{equation}
of the same function (\ref{J6a})
keeps both of its minima and splits the $D_4^-$ point into two saddlepoints and two imaginary critical points, see Fig.~\ref{J14} (center).

The perturbation 
\begin{equation}
x (x^2 + y^4 + \varepsilon y^2 - \varepsilon^3) \label{J4a}
\end{equation}
of singularity $x^3 + x y^4$ of class $J_{10}^1$ has four real critical points: one maximum, one minimum and two saddles, see Fig.~\ref{J14} (left). 

The perturbation 
\begin{equation}
\label{J2a} x^3 + y^6 + \varepsilon (y^2-x)
\end{equation}
has only two real critical points: a minimum and a saddlepoint, see Fig.~\ref{J14} (right).

\begin{figure}
\unitlength 0.6mm
\begin{picture}(40,26)
\put(0,13){\line(1,0){40}}
\put(20,13){\circle{20}}
\put(18,7){\footnotesize $+$}
\put(18,15){\footnotesize $-$}
\put(9.7,13){\circle*{1.33}}
\put(30.3,13){\circle*{1.33}}
\end{picture} \qquad \qquad \quad
\begin{picture}(50,26)
\put(29.8,15){\circle*{1.33}}
\put(20.2,15){\circle*{1.33}}
\bezier{100}(1,4)(10,6)(15,10)
\bezier{300}(15,10)(31,24)(15,26)
\bezier{100}(13,20)(10,25)(15,26)
\bezier{400}(37,20)(25,8)(13,20)
\bezier{100}(35,26)(40,25)(37,20)
\bezier{400}(35,26)(19,24)(35,10)
\bezier{100}(49,4)(40,6)(35,10)
\put(33.00,21.00){\small \makebox(0,0)[cc]{$-$}}
\put(17.00,21.00){\small \makebox(0,0)[cc]{$-$}}
\put(22.5,6){\footnotesize $-$}
\end{picture} \qquad \qquad \quad
\unitlength 0.2mm
\begin{picture}(100,105)
\bezier{800}(0,10)(55,30)(60,80)
\bezier{300}(60,80)(60,94)(50,95)
\bezier{300}(50,95)(40,94)(40,80)
\bezier{800}(100,10)(45,30)(40,80)
\put(50,48){\circle*{6}}
\put(44,72){\scriptsize $-$}
\put(44,10){\scriptsize $-$}
\end{picture}
\caption{Perturbations for $\Phi_1$ with four or two critical points}
\label{J14}
\end{figure}
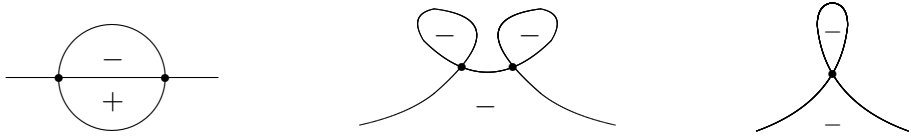

Finally, the perturbation 
\begin{equation}
\label{J0}
 x^3 + y^6 + \varepsilon x
\end{equation} has no real critical points.

\FloatBarrier

\subsection{Realizations of virtual Morse functions of class $\Phi_3$}
Two polynomials of type $\Phi_3$ with ten real critical points, exactly one of which is a maximum point, can be chosen as 
\begin{equation}
x (x+3y^2-2)(x+y^2-1) \label{J10c}
\end{equation}
 and 
\begin{equation}
x (x + 3y(y+1)-1)(x + \frac{1}{3}(y-1)(y-3)-1). \label{J10d}
\end{equation}
They realize the D-graphs shown in Figs.~\ref{77374} and \ref{225148}. 
Their zero-level sets are shown in the left parts of the same figures. The first of these polynomials is invariant under the involution (\ref{invol0}). 

A polynomial of type $\Phi_3$ with ten real critical points, exactly two of which are the maxima, can be chosen as follows: \begin{equation}
x(x+2(y+1)(y-2)-1/2)(x + (y-1)(y+2)-1), \label{J10e}
\end{equation}
see Fig.~\ref{128634}.

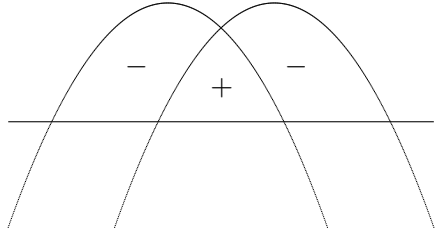
\begin{figure}
\unitlength 0.5 mm
\begin{picture}(80,50)
\put(0,25){\line(1,0){80}}
\bezier{500}(0,5)(30,90)(60,5)
\bezier{500}(20,5)(50,90)(80,5)
\put(22,34){\footnotesize $-$}
\put(52,34){\footnotesize $-$}
\put(38,30){\footnotesize $+$}
\end{picture} 
\caption{$\Phi_3$, eight critical points, one maximum, 66906}
\label{66906}
\end{figure}

A polynomial of type $\Phi_3$ with eight real critical points, exactly one of which is a maximum, is given by 
\begin{equation}
x(x+(y-1)^2-2)(x+(y+1)^2-2)+ \varepsilon y^6 \equiv x^3+2x^2y^2-2x^2+x y^4-6x y^2+x +\varepsilon y^6 \ \label{J8c}
\end{equation}
with any sufficiently small $\varepsilon > 0$.
 The zero level set of this polynomial without the last term $\varepsilon y^6$ 
(which makes its principal quasihomogeneous part non-degenerate of class $J_{10}^3$)
is shown in Fig.~\ref{66906}.

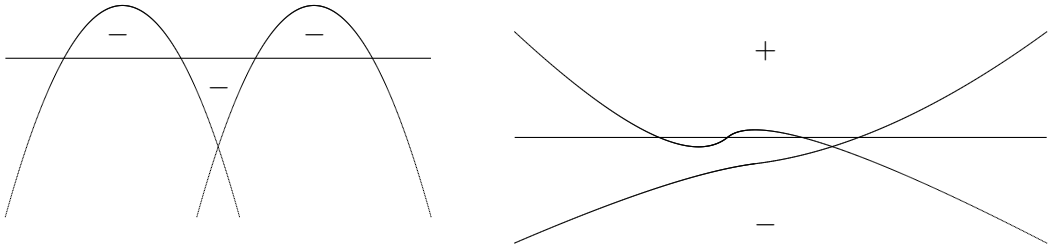
\begin{figure} 
\unitlength 0.55 mm
\begin{picture}(90,50)
\put(0,35){\line(1,0){80}}
\bezier{500}(0,5)(22,85)(44,5)
\bezier{500}(36,5)(58,85)(80,5)
\put(19,38){\footnotesize $-$}
\put(56,38){\footnotesize $-$}
\put(38,28){\footnotesize $-$}
\end{picture} 
 \qquad \qquad
\begin{picture}(100,40)
\put(0,20){\line(1,0){100}}
\bezier{400}(0,40)(30,12)(40,20)
\bezier{400}(40,20)(48,27)(100,0)
\bezier{400}(100,40)(70,18)(45,15)
\bezier{400}(45,15)(30,13)(0,0)
\put(45,35){\footnotesize $+$}
\put(45,2){\footnotesize $-$}
\end{picture}
\caption{$\Phi_3$, eight critical points, no maxima} 
\label{J381}
\end{figure}

Two polynomials representing two different isotopy classes with eight real critical points, none of which are the maxima, are shown in Fig.~\ref{J381}. 
The first of these polynomials is
\begin{equation}
x(x+(y-2)^2-2)(x+(y+2)^2-2)+\varepsilon y^6 \equiv x^3+2x^2y^2+4x^2+x y^4-12xy^2+4x+\varepsilon y^6. \label{J8d}
\end{equation}
It realizes the class with Card = 131148. The Fig.~\ref{J381} (left) shows the zero-level set of this polynomial without the term $\varepsilon y^6$. The second isotopy class has Card = 82350 and no symmetric realizations. Nevertheless, it is achiral, see \S~\ref{achirali}. A realization of this class is shown in Fig.~\ref{J381} (right). To construct it, we first take the perturbation $$x^3 -x y^4 - \varepsilon x^2 y$$ of the original $J_{10}^3$ singularity $x^3 -x y^4$, and then perturb its single real critical point of type $D_8^-$ as shown on p. 15 of \cite{AC}.

An isotopy invariant proving that these two Morse functions indeed belong to different isotopy classes is as follows. For any polynomial $f$ of type $\Phi_3$ with exactly three local minima, a unique parabola (or line) with equation of the form \ $x=a y^2 + by +c$ \ exists that passes through these three points. The restriction of $f$ to this curve is a polynomial of degree at most six in the coordinate $y$. Since it has three minima, it is of degree exactly six. Therefore, the coefficient \ $a$ \ of this parabola is not a root of the polynomial $\alpha t^3 + \beta t^2 + \gamma t + \delta,$ where $\alpha x^3 + \beta x^2 y^2 + \gamma x y^4 + \delta y^6$ is the principal quasihomogeneous part of $f$. Thus, the position of the number \ $a$ \ among three roots of this polynomial is an invariant of the isotopy class. For a polynomial realizing the left picture of Fig.~\ref{J381} the coefficient \ $a$ \ does not separate these roots, while for the right picture it separates one root from the other two.
\medskip

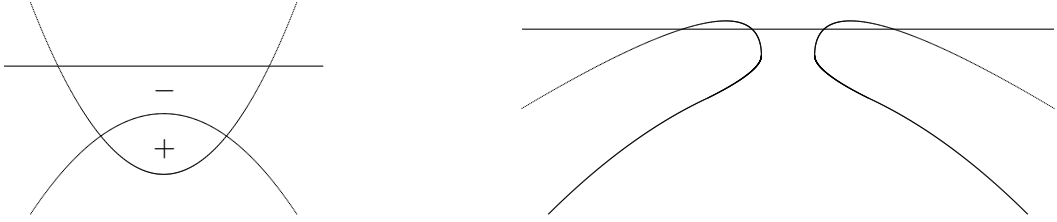
\begin{figure}
\unitlength 0.7mm
\begin{picture}(60,43)
\put(0,28){\line(1,0){60}}
\bezier{400}(5,40)(30,-25)(55,40)
\bezier{300}(5,0)(30,38)(55,0)
\put(28,22){\footnotesize $-$}
\put(28,11){\footnotesize $+$}
\end{picture} \qquad \qquad \qquad
\begin{picture}(100,40)
\put(0,35){\line(1,0){100}}
\put(48,40){\footnotesize $+$}
\put(48,10){\footnotesize $-$}
\bezier{300}(55,30)(55,47)(100,20)
\bezier{300}(55,30)(55,27)(65,22)
\bezier{300}(65,22)(80,15)(95,0)
\bezier{300}(45,30)(45,47)(0,20)
\bezier{300}(45,30)(45,27)(35,22)
\bezier{300}(35,22)(20,15)(5,0)
\end{picture} 
\caption{$\Phi_3$, six critical points; one maximum (left) and no maxima (right)}
\label{J36}
\end{figure}

The isotopy component of type $\Phi_3$ with exactly six real critical points, exactly one (respectively, none) of which is a local maximum, is represented by the polynomial 
\begin{equation}
x(x-y^2+4)(x+y^2+1) \label{J6c}
\end{equation}
(respectively, 
\begin{equation}
(x+\varepsilon)\left(2x^2 +5x y^2+2y^4 -8\frac{\sqrt{2}}{3} x - 10 \frac{\sqrt{2}}{3}y^2 + \frac{25}{9}\right) \label{J6d}
\end{equation}
with sufficiently small $\varepsilon >0$). 
The zero-level sets of these polynomials are shown in Fig.~\ref{J36}.
The degenerate version of the second of these polynomials corresponding to $\varepsilon = 0$ has only two real critical points of type $A_3$. 

\begin{figure}
\unitlength 0.7mm
\begin{picture}(80,43)
\put(0,30){\line(1,0){40}}
\bezier{300}(0,40)(20,0)(40,40)
\bezier{200}(0,0)(20,20)(40,0)
\put(18,23){\footnotesize $-$}
\put(18,13){\footnotesize $+$}
\end{picture} \qquad \qquad \qquad 
\begin{picture}(60,40)
\put(0,20){\line(1,0){40}}
\bezier{300}(0,40)(20,20)(40,40)
\bezier{200}(0,0)(20,20)(40,0)
\put(18,23){\footnotesize $-$}
\put(18,13){\footnotesize $+$}
\end{picture}
\caption{$\Phi_3$. Four critical points, Card=21410 (left); two critical points, Card=14778 (right)}
\label{21410}
\end{figure}
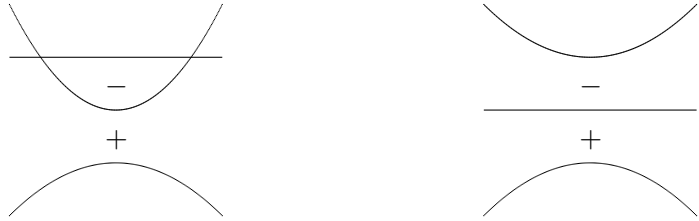

A polynomial of type $\Phi_3$ with exactly four real critical points and no local maxima is given by
\begin{equation}
x (x - y^2+1)(x+y^2 +3),
\label{J4c}
\end{equation}
see Fig.~\ref{21410} (left). A polynomial with only two real critical points 
is given by 
\begin{equation} x (x-y^2-1)(x + y^2 + 1) , \label{J2} 
\end{equation}
see Fig.~\ref{21410} (right).

\section{Chirality of virtual Morse functions. Proof of Theorem \ref{cher}}
\label{chirali}

\subsection{Chirality cocycle and proof of achirality statements of Theorem \ref{cher}.}
\label{achirali}

Denote by $\HH$ the one-dimensional simplicial cochain of the formal graph of type $\Phi_1$ or $\Phi_3$ with coefficients in $\Z_2$, which takes the non-zero value only on the edges corresponding to the type $s2$ elementary operation of changing the order of two neighboring real critical values, {\em both of which are achieved at minimum points}.

\begin{theorem}
\label{caca}
If all virtual Morse functions in a virtual component of $\Phi_1$ or $\Phi_3$ type
 have exactly two or three minima, then this component is chiral if and only if the cocycle $\HH$ is trivial on it.
\end{theorem} 

\noindent
{\it Proof.} Let $N$ be the number of minima of these virtual functions.
Suppose that the virtual component is achiral and hence is associated with only one isotopy class of Morse polynomials in $\Phi_1$ or $\Phi_3$. Then, there exists
a generic path within this isotopy class that connects two mutually symmetric generic polynomials, $f(x,y)$ and $f(x,-y)$. The endpoints of this path are associated with the same virtual Morse function. Therefore, the associated path in the formal graph is a cycle.
Polynomials of classes $\Phi_1$ and $\Phi_3$ cannot have two minima on the same vertical line $\{x=c\}$ because they are of degree three in the $x$ variable.
Thus, our path preserves the order of the projections of the minimum points to the $y$ axis. Consider the following two permutations of the numbers $1, 2, \dots, N$. The permutation $\sigma$ depends only on the initial polynomial $f$ and maps each number $i$ to the order of the critical value of $f$ at the $i$-th minimum point from the left among all such critical values. The permutation $\varkappa$ is the permutation of orders of critical values at minima defined by the continuation along our path. The elementary transpositions that form the permutation $\varkappa$ along this path correspond to the edges of the associated path in the virtual component, on which the cochain $\HH$ takes the non-zero value.

For any $i \in [1, N]$ the critical value of the $i$th minimum point from the left of the initial polynomial of this path is equal to the critical value of the $(N-i)$-th minimum point of the final polynomial. Therefore, the permutations $\varkappa$ and 
{\tiny $\left(\begin{array}{ccccc}
1& 2 & \dots & N-1 & N\\
N & N-1 & \dots & 2 & 1
\end{array} \right)$}
 are conjugate via $\sigma$, in particular, have the same parity. For $N=2$ or 3, this parity is odd, and hence the cocycle $\HH$ takes a non-zero value on the cycle. 

Conversely, given a cycle in a virtual component of the formal graph, we can realize it by a path in the space $\Phi_1$ or $\Phi_3$ of real polynomials starting from an arbitrary vertex of the cycle. That is, the path begins with a polynomial associated with this vertex and successively undergoes all surgeries encoded by the edges of the cycle. The existence of such a path is ensured by Proposition \ref{propmain}. If the cocycle $\HH$ takes the non-zero value on the initial cycle, then the isotopy classes of the endpoints of the realized path
in the space of Morse polynomials differ because a closed loop in $\Phi_1$ or $\Phi_3$ defines the trivial permutation $\varkappa$. 
\hfill $\Box$

\begin{remark} \rm We could formulate the same criterion in terms of the critical values at the maximum points. Therefore, if both inertia indices $m_-$ and $m_+$ of the passport invariant of the virtual component are equal to two or three, then another cocycle defined in this virtual component, similar to $\HH$ but in the terms of maxima instead of minima, is homologous to $\HH$.
\end{remark}

Almost all of the virtual components that Theorem \ref{cher} claims are achiral are represented by symmetric polynomials with respect to the involution (\ref{invol0}), see the polynomials (\ref{J10b})--(\ref{J8b}), (\ref{J6a})--(\ref{J10c}), and (\ref{J8c})--(\ref{J2}). The remaining two virtual Morse functions are of type $\Phi_3 $ with eight real critical points and Card = 82350. One of them has no maxima (see Fig.~\ref{J381} right) and the other is obtained from it by the involution (\ref{invol}). The simplicial cocycle $\HH$ takes the non-zero value on the cycle $A \to B \to D \to C \to A$ in the corresponding virtual component, see Fig.~\ref{cocy} and Example \ref{cocycleh}. Indeed, this cocycle takes the non-zero value only on the edge $[A, C]$ of this cycle. This proves the achirality of these two virtual components as well.

\begin{remark} \rm
This cycle was found by a minor suspension of our program.
\end{remark}

\begin{remark} \rm Another virtual component with eight real critical points, no maxima and Card=131148  contains the virtual Morse function 
\begin{equation}
\label{altar}
\begin{array}{|ccc|ccccc|cc|}
\hline 
$-2$  &  0  &  0  &  0  &  0  &  1  &  1  &  0  & 0 & 0 \\
  0  &  $-2$  &  0  &  1 &   0  &  0  &  1  &  1 &  $1$ &  $ 1$ \\
  0  & 0 &  $-2$ &   1 &   1  &  0  &  0  &  0  &  1 &  1 \\
  0  &  1 &   1 &  $-2$ &   0  &  0  &  0  &  0  &  $-1$  &  $-1$ \\
  0   & 0  &  1  &  0  & $-2$ &  0  &  0  &  0  &  0  &  0 \\
  1  &  0  &  0  &  0   & 0  & $-2$  &  0  &  0  &  0  &  0 \\
  1  &  1  &  0  &  0  &  0  &  0 &  $-2$  &  0  & $-1$  &  $-1$ \\
  0  &  1  &  0  &  0  &  0  &  0  &  0 &  $-2$ &  0  &  0 \\
 $0$ &  $1$ &  1  & $-1$ &   0  &  0  &  $-1$  &  0 &  $-2$ & $-2$ \\
 $ 0$ & $1$ &  1  & $-1$  &  0   & 0  &  $-1$  &  0  & $-2$ &  $-2$ \\
\hline
    2 &   2  &  2  & $-2$ & $-2$ & $-2$ &  $-2$ &  $-2$ &  $-2$  & $-2$ \\
\hline
    3  &  3  &  3  &  2   & 2  &  2  &  2 &   2   & & \\
\hline
\end{array}
\end{equation}
In this virtual component, a cycle that  proves its achirality consists of only $s2$ type edges.  It defines the  transpositions of critical values at the critical points 1 and 3, 5 and 6, and 4 and 7. See also Fig.~\ref{J381} (left).
\end{remark}

\subsection{Proof of chirality statements of Theorem~\ref{cher}}
\label{cchirali}

\begin{proposition}
If a virtual component of type $\Phi_1$ or $\Phi_3$ consists of virtual functions with ten real critical points and is achiral, then its $D$-graph has an automorphism, whose restriction to the set of vertices corresponding to minimum points $($or to the set of vertices corresponding to maximum points$)$ is an involution with at most one fixed point.
\end{proposition}

\noindent
{\it Proof.} Let $f(x, y)$ be a generic Morse polynomial in an achiral isotopy class with ten real critical points. There are two one-to-one correspondences between its critical points and the critical points of the polynomial $f(x, -y)$. One correspondence preserves the critical values of the critical points, and the other is obtained by tracing the critical points along the path connecting the two polynomials in the isotopy class. The composition of these correspondences of critical points (and, consequently, of the corresponding vertices of the associated D-graphs) extends to an automorphism of the D-graph of $f$. This automorphism maps the vertex corresponding to $i$-th from the right 
minimum point to the vertex corresponding to the $i$-th from the left minimum point. Similarly, it permutes the maximum points. In particular, it defines an involutive permutation of each of these two sets of points with at most one fixed point in each. \hfill $\Box$
\medskip

This proposition implies the chirality of all virtual components with ten real critical points, listed in Theorems \ref{enu1} and \ref{enu3}, except for those whose achirality follows from formulas (\ref{J10b}) and (\ref{J10c}), see Figs.~\ref{26378} and \ref{77374}. 

It remains to prove the chirality of the isotopy class of type $\Phi_1$ with eight real critical points, none of which are maxima, see Fig.~\ref{J38} (left). (The criterion of Theorem \ref{caca} does not work for this class because the number of minima is not 2 or 3). 

An isotopy invariant that separates each polynomial $f(x,y) $ of this class from its mirror image $f(x, -y)$ is as follows. There is a single polynomial equation \ $x = a y^3+b y^2+c y +d$ \ whose graph in $\R^2$ contains the  four minimum points of $f$. The coefficient \ $a$ \ of this polynomial is never zero: otherwise, the restriction of the function $f$ to this  curve would be a degree six polynomial in $y$ with four minimum points. Therefore, the sign of this coefficient is an invariant of the isotopy class. These coefficients defined by the functions $f(x,y)$ and $f(x, -y)$ are opposite.
 \hfill $\Box$ \medskip

 This concludes the proof of Theorem \ref{cher}. \hfill $\Box$

\section{Proof of Theorem \ref{tabadjp}}
\label{sectlast} 

Recall that the {\em canonical Coxeter-Dynkin graphs of the real singularity classes $D_4^+$, $D_6^+$, and $D_8^+$} are the graphs given in Fig.~\ref{ddr}.
\begin{figure}
\unitlength 0.43mm
\begin{picture}(20,28)
\put(0,15){\circle*{1.5}}
\put(0,15){\line(1,1){12}}
\put(0,15){\line(1,-1){12}}
\put(24,15){\circle*{1.5}}
\put(24,15){\line(-1,-1){12}}
\put(24,15){\line(-1,1){12}}
\put(2, 15){\line(1,0){5}}
\put(9,15){\line(1,0){5}}
\put(16,15){\line(1,0){5}}
\put(12,27){\circle*{1.5}}
\put(12,3){\circle*{1.5}}
\end{picture} \quad
\begin{picture}(110,30)
\put(30,15){\line(1,0){30}}
\put(30,15){\circle*{1.5}}
\put(45,15){\circle*{1.5}}
\put(60,15){\circle*{1.5}}
\put(60,15){\line(1,1){12}}
\put(60,15){\line(1,-1){12}}
\put(84,15){\circle*{1.5}}
\put(84,15){\line(-1,-1){12}}
\put(84,15){\line(-1,1){12}}
\put(62, 15){\line(1,0){5}}
\put(69,15){\line(1,0){5}}
\put(76,15){\line(1,0){5}}
\put(72,27){\circle*{1.5}}
\put(72,3){\circle*{1.5}}
\end{picture} \quad
\begin{picture}(100,30)
\put(0,15){\line(1,0){60}}
\put(0,15){\circle*{1.5}}
\put(15,15){\circle*{1.5}}
\put(30,15){\circle*{1.5}}
\put(45,15){\circle*{1.5}}
\put(60,15){\circle*{1.5}}
\put(60,15){\line(1,1){12}}
\put(60,15){\line(1,-1){12}}
\put(84,15){\circle*{1.5}}
\put(84,15){\line(-1,-1){12}}
\put(84,15){\line(-1,1){12}}
\put(62, 15){\line(1,0){5}}
\put(69,15){\line(1,0){5}}
\put(76,15){\line(1,0){5}}
\put(72,27){\circle*{1.5}}
\put(72,3){\circle*{1.5}}
\end{picture}
\caption{Real Coxeter-Dynkin graphs $D_4^+$, $D_6^+$, and $D_{8}^+$}
\label{ddr}
\end{figure}
The canonical Coxeter-Dynkin graph of the real singularity $D_{2k}^-$ is just the standard $D_{2k}$ graph. The canonical Coxeter-Dynkin graphs of the other real simple singularities, $A_k$, $D_{2k+1}$, $E_6$, $E_7$, and $E_8$, are also the same as in the usual ``complex'' theory. 

All critical points of small perturbations of the $J_{10}$ singularities are simple or belong to the $J_{10}$ class itself.

\begin{proposition}
\label{pro18}
The following two conditions are equivalent:

\begin{enumerate}
\item
there exists a polynomial of type $\Phi_1$ $($respectively, $\Phi_3)$ having two real critical points of classes $\Xi$ and $\tilde \Xi$ with $\mu(\Xi)+\mu(\tilde \Xi)=10$;

\item
the set of vertices of one of the D-graphs shown in Figs.~\ref{122298}--\ref{33528} $($respectively, \ref{77374}--\ref{29370}$)$ can be divided into two subsets of cardinality \ $\mu(\Xi)$ \ and \ $\mu(\tilde \Xi)$ \ in such a way that 
\begin{enumerate}
\item
all edges of the D-graph connecting vertices of different subsets are directed from vertices of the first subset to vertices of the second, and
\item
the edges, whose vertices both belong to one of these subsets, form the canonical Coxeter-Dynkin graph of the corresponding real simple singularity $\Xi$ or $\tilde \Xi$. 
\end{enumerate}
\end{enumerate}
\end{proposition}

\noindent
{\it Proof.} A proof of the ``only if'' part repeats the proof of Proposition 26 in \cite{Vx9}. The ``if'' part follows from the same considerations regarding the Lyashko--Looijenga covering, cf. \cite{Lya}, \cite{Jaw2}. Namely, let $f$ be a polynomial of class $\Phi_1$ or $\Phi_3$ that realizes this D-graph. For certainty, suppose that the edges of the D-graph are directed from the vertices of the $\Xi$-subgraph to the vertices of the $\tilde \Xi$-subgraph.
Using the group $\G$, we can assume that $f$ has the form (\ref{vers1}) or (\ref{vers3}). Using the Lyashko--Looijenga map as in \cite{Vx9}, we can continuously deform this polynomial into a generic polynomial $\tilde f$, such that all of its critical values corresponding to vertices of the $\Xi$-subgraph lie below all critical values of the $\tilde \Xi$-subgraph. (This deformation may intersect the set of non-generic polynomials, at which the critical values of distant critical points are equal; however, this does not prevent the deformation.) Let $(c_1<c_2< \dots <c_{10}) \in \R^{10}$ be the set of critical values of $\tilde f$. Consider the parametrized segment $p: [0,1] \to \R^{10}$ that connects this point with the point $(c_-, \dots, c_-, c_+, \dots, c_+)$ consisting of $\mu(\Xi)$ copies of the mean value $c_-$ of numbers $c_1, \dots, c_{\mu(\Xi)}$ and $\mu(\tilde \Xi)$ copies of the mean value $c_+$ of numbers $c_{\mu(\Xi)+1}, \dots, c_{10}$. According to Proposition 2 of \cite{Jaw2}, this segment can be lifted to a path $ \{\tilde f_t\},$ $t\in [0,1]$, in the space of polynomials (\ref{vers1}) or (\ref{vers3}) such that $\tilde f_0\equiv \tilde f$ and  $p(t)$ is the collection of critical values of the polynomial $\tilde f_t$ for any $t\in [0,1]$. The final point $\tilde f_1$ of this path is the desired polynomial with the critical points of classes $\Xi$ and $\tilde \Xi$. \hfill $\Box$ \medskip

The condition 1) of Proposition \ref{pro18} is equivalent to the condition $\{\Xi + \tilde \Xi\}\rightsquigarrow J_{10}^1$ or  $\{\Xi + \tilde \Xi\}\rightsquigarrow J_{10}^3)$ from Theorem \ref{tabadjp} because the spaces $\Phi_1$ and $\Phi_3$ are versal deformations of $J_{10}^1$ and $J_{10}^3$ singularities, and the group (\ref{dila}) provides the functions with desired critical points arbitrarily close to the origin. Therefore, Theorem \ref{tabadjp} follows immediately from Proposition \ref{pro18} and the lists of D-graphs given in Theorems \ref{enu1} and \ref{enu3}. \hfill $\Box$

\section{A problem}

To explicitly present the polynomials representing all isotopy classes with ten real critical points. This has only been done above for classes whose $D$-graphs are shown in Figs.~\ref{122298}, \ref{26378}, \ref{3852} (left), \ref{77374}, \ref{225148}, \ref{128634}, and \ref{59862} (right).

}

\begin{thebibliography}{99}

\bibitem{AC} N.~A'Campo, {\it Le groupe de monodromie du déploiement des singularités isolées de courbes planes. I,} Math. Ann. 213, 1-32 (1975), doi: 10.1007/BF01883883

\bibitem{Acongr} V.I.Arnold, {\it Critical points of smooth functions}, Proc. of the International Congress of Mathematicians (Vancouver 1974), 19--40.

\bibitem{Kluwer} V.I.~Arnold, {\it Singularities of caustics and wave fronts}. Vol. 62. Springer Science \& Business Media, 2001.

\bibitem{AVGZ} V.I.~Arnold, S.M.~Gusein--Zade, A.N.~Varchenko, {\it Singularities of differentiable maps}, Vols. 1 and 2, Birkhäuser, Basel, 2012.

\bibitem{AGLV2} V.I.~Arnol’d, V.V.~Goryunov, O.V.~Lyashko, V.A.~Vassiliev, 
{\it Singularity Theory. II: Classification and Applications}. VINITI, 1989, 5–249. Engl. Transl.: Encyclopaedia of Mathematical Sciences. 39. Berlin: Springer-Verlag, 1993, 235 p.

\bibitem{GZ} S.M.~Gusein--Zade, {\it Intersection matrices for some singularities of functions of two variables}, Funct. Anal. Appl. 8:1 (1974), 10-13.

\bibitem{Jaw} P.~Jaworski, {\it Distribution of critical values of miniversal deformations of parabolic singularities}, Invent. Math., 1986, 86:1, 19--33.

\bibitem{Jaw2} P.~Jaworski, {\it Decompositions of parabolic singularities}, Bull. Sci. Math. (2) 112:2 (1988), 143–176.

\bibitem{Liv} I.S.~Livshits, Automorphisms of the complement of the bifurcation set of functions for simple singularities, Funct. Anal. Appl. 15:1 (1981), 29–32.

\bibitem{Lya} O.V.~Lyashko, {\it Decompositions of simple singularities of functions}. Funct. Anal. Appl. 10:2 (1976), 122-127.

\bibitem{Lo0} E.~Looijenga, {\it The complement of the bifurcation variety of a simple singularity}, Invent. Math. 23 (2), 105--116.

\bibitem{M} J.~Milnor, {\it Singular points of complex hypersurfaces,} Princeton University Press (1968).

\bibitem{sed} V.D.~Sedykh, {\it On the topology of stable Lagrangian maps with singularities of types $A$ and $D$}, Izvestiya: Mathematics,
2015, Volume 79, Issue 3, 581–622.

\bibitem{sede} V.D.~Sedykh, {\it The topology of the complement to the caustic of a Lagrangian germ of type $E_6^\pm$}, Russian Math. Surveys, {\bf 78:3} (2023), 569--571.

\bibitem{thom} R.~Thom, {\it Topological models in biology}, Topology, 8: 3 (1969), 313--335.

\bibitem{thom2} R.~Thom, {\it The bifurcation subset of a space of maps}, Manifolds, Amsterdam, 1970, Springer Lecture Notes in Math. 197 (1971), 202—208.

\bibitem{APLT} V.A.~Vassiliev, {\it Applied Picard-Lefschetz theory,} AMS, Providence RI, 2002.

\bibitem{VS} V.A.~Vassiliev, {\it Real Function Singularities and Their Bifurcation Sets}, in: Handbook of Geometry and Topology of Singularities VII, eds. J.L.~Cisneros-Molina, L\^e D\~ung Tr\'ang, J.~Seade, Springer, 2025, 71–119

\bibitem{Vsing} V.A.~Vassiliev, {\it Complements of caustics of real function singularities}, Journal of Singularities, 27 (2024), 47–67 , arXiv: 2304.09824 

\bibitem{Vx9} V.A.~Vassiliev, {\it Isotopy classification of Morse polynomials of degree four in ${\mathbb R}^2$}, Moscow Math. Journal, 25:2 (2025), 249–299 arXiv: 2311.11113

\bibitem{Vp8} V.A.~Vassiliev, {\it Isotopy classification of Morse polynomials of degree 3 in ${\mathbb R}^3$}, 2024, arXiv: 2404.17891





\end{thebibliography}
\end{document}